\documentclass[reqno]{amsart}

\usepackage[utf8]{inputenc}
\usepackage[english]{babel} 
\usepackage{amsmath,amssymb,amsfonts,amsthm}
\usepackage{mathrsfs, esint, dsfont}
\usepackage{moreverb,rotating,graphics,float}
\usepackage[colorlinks, linkcolor={blue!50!blue}, citecolor={red}]{hyperref}
\usepackage{framed,fancybox}
\usepackage{enumitem}
\usepackage{csquotes}
\usepackage{todonotes}

\newtheorem{theorem}{Theorem}
\newtheorem{proposition}[theorem]{Proposition}
\newtheorem{remark}[theorem]{Remark}
\newtheorem{lemma}[theorem]{Lemma}
\newtheorem{corollary}[theorem]{Corollary}

\newcommand\1{{\mathds 1}}

\def\bP{{\mathbb P}}
\def\bbI{{\mathbb I}}
\def\N{{\mathbb N}}

\def\R{{\mathbb R}}

\def\rd{{\mathrm{d}}}
\def\re{{\mathrm{e}}}

\def\cB{{\mathcal B}}

\def\cE{{\mathcal E}}

\def\gS{{\mathfrak{S}}}
\newcommand{\ii}{\infty}
\newcommand{\dps}{\displaystyle}
\newcommand{\nn}{\nonumber}
\newcommand{\Tr}{{\rm Tr}}

\newcommand{\tr}{{\rm Tr}}
\newcommand{\bra}{\langle}
\newcommand{\ket}{\rangle}

\renewcommand{\epsilon}{\varepsilon}

\newcommand\pscal[1]{{\ensuremath{\left\langle #1 \right\rangle}}}
\newcommand{\norm}[1]{ \left| \! \left| #1 \right| \! \right| }

\DeclareMathOperator{\rank}{Rank}
\newcommand{\wto}{\rightharpoonup}

\author{Rupert L. Frank}
\address[Rupert L. Frank]{Mathematics 253-37, Caltech, Pasa\-de\-na, CA 91125, USA, and Mathematisches Institut, Ludwig-Maximilans Universit\"at M\"unchen, Theresienstr. 39, 80333 M\"unchen, Germany, and Munich Center for Quantum Science and Technology (MCQST), Schellingstr. 4, 80799 M\"unchen, Germany}
\email{rlfrank@caltech.edu}

\author{David Gontier}
\address[David Gontier]{CEREMADE, University of Paris-Dauphine, PSL University, 75016 Paris, France}
\email{gontier@ceremade.dauphine.fr}

\author{Mathieu Lewin}
\address[Mathieu Lewin]{CNRS and CEREMADE, University of Paris-Dauphine, PSL University, 75016 Paris, France}
\email{mathieu.lewin@math.cnrs.fr}

\title[NLS equation for orthonormal functions II]{The nonlinear Schr\"odinger equation for orthonormal functions\\ II. Application to Lieb-Thirring inequalities}

\date{August 27, 2020}

\begin{document}

\begin{abstract}
In this paper we disprove part of a conjecture of Lieb and Thirring concerning the best constant in their eponymous inequality. We prove that the best Lieb-Thirring constant when the eigenvalues of a Schr\"odinger operator $-\Delta+V(x)$ are raised to the power $\kappa$ is never given by the one-bound state case when $\kappa>\max(0,2-d/2)$ in space dimension $d\geq1$. When in addition $\kappa\geq1$ we prove that this best constant is never attained for a potential having finitely many eigenvalues. The method to obtain the first result is to carefully compute the exponentially small interaction between two Gagliardo-Nirenberg optimisers placed far away. For the second result, we study the dual version of the Lieb-Thirring inequality, in the same spirit as in Part I of this work~\cite{GonLewNaz-20_ppt}. In a different but related direction, we also show that the cubic nonlinear Schr\"odinger equation admits no orthonormal ground state in 1D, for more than one function. 

\bigskip

\noindent \sl \copyright~2020 by the authors. This paper may be reproduced, in its entirety, for non-commercial purposes.
\end{abstract}

\maketitle

\section{Introduction and main results}

This paper is a continuation of a previous work~\cite{GonLewNaz-20_ppt} where the last two authors together with F.Q.~Nazar studied the existence of ground states for the \emph{nonlinear Schr\"odinger equation} (NLS) for systems of orthonormal functions. In the present paper, we exhibit a connection between the corresponding minimisation problem and the family of Lieb-Thirring inequalities~\cite{LieThi-75,LieThi-76,LieSei-09}, which enables us to prove results both for the Lieb-Thirring inequalities and the NLS equation studied in~\cite{GonLewNaz-20_ppt}. 

\subsection{Lieb-Thirring inequalities}\label{sec:LT}
The Lieb-Thirring inequality~\cite{LieThi-75,LieThi-76} is one of the most important inequalities in mathematical physics. It has been used by Lieb and Thirring~\cite{LieThi-75} 
to give a short proof of the stability of matter~\cite{DysLen-67,DysLen-68,Lieb-90,LieSei-09} and it is a fundamental tool for studying large fermionic systems. It is also a source of many interesting mathematical questions.  

\subsubsection{The finite rank Lieb-Thirring constant}
Let $d\geq 1$, $\kappa \geq 0$ and $N\geq 1$, and let  $L_{\kappa,d}^{(N)}$ be the best constant in the {\em finite rank Lieb-Thirring inequality}
\begin{equation}
    \boxed{ \sum_{n=1}^N |\lambda_n(-\Delta+V)|^\kappa \leq L_{\kappa,d}^{(N)} \int_{\R^d} V(x)_-^{\kappa+\frac{d}2}\,{\rd x} }
\label{eq:LT_V_N}
\end{equation}
for all $V\in L^{\kappa+\frac{d}{2}}(\R^d)$, where $a_-=\max(0,-a)$ and $\lambda_n(-\Delta+V)\leq0$ denotes the $n$th min-max level of $-\Delta+V$ in $L^2(\R^d)$, which equals the $n$th negative eigenvalue (counted with multiplicity) when it exists and 0 otherwise. The constant $L_{\kappa, d}^{(1)}$ is finite by the Gagliardo-Nirenberg inequality, under the assumption that 
\begin{equation}
\begin{cases}\kappa\geq\frac12& \text{in $d=1$,}\\ 
\kappa>0& \text{in $d=2$,}\\
\kappa\geq0& \text{in $d\geq3$.}
\end{cases}
\label{eq:constraint_kappa}
\end{equation}
These restrictions on $\kappa$ are optimal in the sense that $L_{\kappa,d}^{(1)}=\infty$ for $0\leq\kappa<1/2$ in $d=1$ and for $\kappa=0$ in $d=2$. Note that $L_{\kappa, d}^{(N)}$ is finite under the same restrictions as for $L_{\kappa, d}^{(1)}$, since $L_{\kappa, d}^{(N)}\leq NL_{\kappa, d}^{(1)}$. Moreover, from the definition we have $L_{\kappa, d}^{(N)} \le L_{\kappa, d}^{(N+1)}$. The Lieb-Thirring theorem states that the limit is finite:
\begin{equation}
L_{\kappa,d} :=L_{\kappa,d} ^{(\ii)}= \lim_{N\to\ii}L_{\kappa,d}^{(N)}<\ii\qquad\text{for $\kappa$ as in~\eqref{eq:constraint_kappa}.}
\label{eq:LT_V_form}
\end{equation}
This was proved by Lieb and Thirring~\cite{LieThi-75,LieThi-76} for $\kappa>1/2$ in $d=1$ and for $\kappa>0$ in $d\geq 2$. The critical cases
$\kappa=0$ in $d\geq3$ and $\kappa=1/2$ in $d=1$ are respectively due to  Cwikel-Lieb-Rozenblum~\cite{Cwikel-77,Lieb-76b,Rozenbljum-72} and Weidl~\cite{Weidl-96}. 

An important question is to determine the value of the optimal Lieb-Thirring constant $L_{\kappa,d}$. This plays for instance a central role in Density Functional Theory~\cite{LewLieSei-19_ppt}. One possibility is that it is attained for an optimal potential $V$ having $N<\ii$ bound states, that is, $L_{\kappa,d}=L^{(N)}_{\kappa,d}$. An opposite scenario is that a sequence $V_N$ of optimal potentials for $L_{\kappa,d}^{(N)}$ tends to be very spread out and flat as $N\to\ii$ so as to have more and more bound states. In this case $L_{\kappa,d}$ is equal to the \emph{semi-classical constant}
\begin{equation} \label{eq:L_sc}
    L_{\kappa, d}^{\rm sc} := \frac{\Gamma\left(\kappa+1\right)}{2^d\pi^{\frac{d}2}\,\Gamma\left(\kappa+d/2+1\right)}.
\end{equation}
Indeed, recall that if we scale a fixed nice potential $V$ with $V_-\neq0$ in the manner $V(\hbar x)$, we obtain in the limit $\hbar\to0$
\begin{align*}
\lim_{\hbar\to0}\frac{\sum_{n\geq1}\big|\lambda_n\big(-\Delta+V(\hbar\cdot)\big)\big|^\kappa}{\int_{\R^d}V(\hbar x)^{\kappa+d/2}_-\,\rd x}&=\lim_{\hbar\to0}\frac{\hbar^d\sum_{n\geq1}\big|\lambda_n(-\hbar^2\Delta+V)\big|^\kappa}{\int_{\R^d}V(x)^{\kappa+d/2}_-\,\rd x}\\
&=\frac{\iint_{\R^{2d}}\big(|p|^2+V(x)\big)_-^\kappa\rd x\,\rd p}{(2\pi)^{d}\int_{\R^d}V(x)^{\kappa+d/2}_-\,\rd x}\\
&=(2\pi)^{-d}\int_{\R^d}(|p|^2-1)_-^\kappa\rd p =L_{\kappa,d}^{\rm sc}.
\end{align*}
Lieb and Thirring conjectured in~\cite{LieThi-76} that the best constant should be given either by the one bound state  case, or by semi-classical analysis:
\begin{equation}
 L_{\kappa,d}\overset{?}{=}\max\left\{L_{\kappa,d}^{(1)},L_{\kappa,d}^{\rm sc}\right\}.
 \label{eq:Lieb-Thirring_conjecture}
\end{equation}
This conjecture has generated a huge interest in mathematical physics. Although the conjecture is still believed to hold in dimension $d=1$, it is now understood that the situation is more complicated in dimensions $d\geq2$. In Section~\ref{sec:discussion} below we will give a precise account of what is known and what is not as of today. Most of the previous works have focused on determining when $L_{\kappa,d}$ equals the semi-classical constant $L_{\kappa,d}^{\rm sc}$. Much fewer works have studied the plausibility that $L_{\kappa,d}$ equals $L^{(1)}_{\kappa,d}$ or even $L^{(N)}_{\kappa,d}$ for some $N\geq1$. In the next section we state our results in this direction.

\subsubsection{Results on the non-optimality of the finite rank Lieb-Thirring constant}
Our first theorem states that for an appropriate range of $\kappa$, the optimal constant in the Lieb-Thirring inequality can never be attained by a potential having finitely many bound states.

\begin{theorem}[Non optimality of the finite-rank case]\label{thm:LT}
    Let $d\geq1$ and
    \begin{equation}
     \begin{cases}
    \kappa >\frac32&\text{for $d=1$,}\\
    \kappa >1&\text{for $d=2$,}\\
    \kappa \geq1&\text{for $d\geq3$.}
    \end{cases}
    \label{eq:condition_kappa}
    \end{equation}
    Then there exists an infinite sequence of integers 
    $N_1 = 1 < N_2 = 2 <N_3<\cdots$
    such that 
    \begin{equation*}
    L^{(N_k-1)}_{\kappa,d}<L^{(N_{k})}_{\kappa,d} \qquad\text{for all}\ k\geq 1.
    \end{equation*}
    In particular, we have 
    $$\boxed{L_{\kappa,d}^{(N)} < L_{\kappa,d}\qquad\text{for all $N\geq1$.}}$$
    In addition, for any $N\geq2$ there exist optimisers $V_N$ for $L_{\kappa,d}^{(N)}$. When $N=N_k$ we have $\lambda_N(-\Delta+V_N)<0$, that is, $-\Delta+V_N$ has at least $N$  negative eigenvalues.
\end{theorem}

As we will discuss below, this result, in particular, disproves the Lieb--Thirring conjecture~\eqref{eq:Lieb-Thirring_conjecture} in dimension $d=2$ in the range $1<\kappa\lesssim 1.165$ and suggests a new scenario for the optimal Lieb-Thirring constant.

It is unclear whether the passage to a subsequence is really necessary or whether the conclusion holds also for $N_k=k$.

The proof of Theorem~\ref{thm:LT} proceeds by studying the \emph{dual formulation} of the Lieb-Thirring inequality~\eqref{eq:LT_V_N} in a similar manner as what was done in~\cite{GonLewNaz-20_ppt} for the nonlinear Schr\"odinger equation. This is explained in detail in the next section, where we also collect more properties of $V_N$. This duality argument requires the assumption $\kappa\geq 1$. It is an interesting open question whether Theorem~\ref{thm:LT} is valid for all 
$\kappa>\max\{0,2-d/2\}$ instead of~\eqref{eq:condition_kappa}. The value of the critical exponent $\max\{0,2-d/2\}$ will be motivated later. In Section~\ref{sec:proof_LT_V_kappa<1} we provide a direct proof for $N=2$ which covers this range of $\kappa$, as stated in the following result. 

\begin{theorem}[Non optimality of the $N=1$ case]\label{thm:LT_bis}
    Let $d\geq1$ and
    \begin{equation}
    \kappa >\max\left\{0,2-\frac{d}2\right\}.
    \label{eq:condition_kappa_bis}
    \end{equation}
    Then we have 
    \begin{equation*}
\boxed{ L^{(1)}_{\kappa,d}<L^{(2)}_{\kappa,d}\leq L_{\kappa,d}.}
    \end{equation*}
\end{theorem}

As we will discuss below, this result, in particular, disproves the Lieb--Thirring conjecture~\eqref{eq:Lieb-Thirring_conjecture} in dimension $d=3$ in the range $1/2<\kappa\lesssim 0.8627$. Thus, together with a result of Helffer-Robert~\cite{HelRob-90} recalled below, the Lieb--Thirring conjecture~\eqref{eq:Lieb-Thirring_conjecture} in dimension $d=3$ is now disproved in the range $1/2<\kappa<1$.

The conclusion $L^{(1)}_{\kappa,d}< L_{\kappa,d}$ for the appropriate range of $\kappa$ is new for all dimensions $2\leq d\leq 7$. Let us briefly sketch an alternative way of arriving at this strict inequality for $d\geq 8$. It is shown in~\cite{GlaGroMar-78} that the best Cwikel-Lieb-Rozenblum constant satisfies $L_{0,d}>L_{0,d}^{\rm sc} > L^{(1)}_{0,d}$ in dimensions $d\geq8$; see also~\cite{Frank-20_ppt}. Here, the constant $L^{(1)}_{0,d}$ is defined in terms of the Sobolev optimiser. A variation of the monotonicity argument from~\cite{AizLie-78} shows that $L^{(1)}_{\kappa,d}/L^{\rm sc}_{\kappa,d}$ is strictly decreasing (see Theorem~\ref{thm:crossing} and Lemma~\ref{al} below). This implies that $L_{\kappa,d}\geq L_{\kappa,d}^{\rm sc}>L^{(1)}_{\kappa,d}$ for all $\kappa\geq0$ and all $d\geq 8$, as claimed. In contrast to this argument, our Theorem~\ref{thm:LT_bis} is not only valid in all dimensions, in the mentioned range of $\kappa$, but it gives the additional information that the two-bound states constant $L_{\kappa,d}^{(2)}$ is above $L_{\kappa,d}^{(1)}$. The mechanism used in our proof is completely different from~\cite{GlaGroMar-78}. There, the authors increased the coupling constant in front of the potential to reach the semi-classical limit. On the other hand, the proof of Theorem~\ref{thm:LT_bis} consists of placing two copies of the one-bound state optimiser far away in the appropriate manner, and computing the resulting exponentially small attraction.

Our proof of Theorem \ref{thm:LT_bis} does not work for $\kappa=0$ in dimensions $d=5,6,7$ (where one still has $2-d/2<0$). Understanding this case is an open problem.

\subsubsection{Discussion}\label{sec:discussion}

We now discuss in detail the consequences of Theorems~\ref{thm:LT} and~\ref{thm:LT_bis} with regard to the Lieb-Thirring conjecture~\eqref{eq:Lieb-Thirring_conjecture}. 

There are many results on the Lieb-Thirring constants $L_{\kappa,d}$. The best estimates currently known are in~\cite{FraHunJexNam-18_ppt}. We mention here a selection of results pertinent to our theorem and refer to~\cite{Frank-20_ppt} for a detailed discussion of known results and open problems. We recall the following known properties:
\begin{itemize}[leftmargin=*]
    \item (Lower bound~\cite{LieThi-76})  For all $d \ge 1$, $\kappa \ge 0$, we have 
    \begin{equation}
L_{\kappa, d} \ge \max \left\{ L_{\kappa, d}^{(1)}, L_{\kappa, d}^{\rm sc} \right\};
\label{eq:LT_conjecture}
    \end{equation}
    \item (Monotonicity~\cite{AizLie-78}) For all $d \ge 1$ and all $1\leq N\leq\ii$, the map $\kappa\mapsto L^{(N)}_{\kappa,d}/L_{\kappa,d}^{\rm sc}$ is non-increasing;\footnote{Only the case $N=\ii$ is considered in~\cite{AizLie-78} but the argument applies the same to any finite $N\geq1$. For $N=1$, we will see in Theorem~\ref{thm:crossing} that $\kappa\mapsto L^{(1)}_{\kappa,d}/L_{\kappa,d}^{\rm sc}$ is in fact \emph{strictly decreasing}.}
    \item ($\kappa=3/2$ in $d=1$~\cite{LieThi-76}) In dimension $d = 1$ with $\kappa = \frac32$, we have, for all $N \in \N$,
    \begin{equation}
    L_{3/2,1}=L_{3/2,1}^{(N)}=L_{3/2,1}^{\rm sc};
    \label{eq:1D_3/2}
    \end{equation}
    \item ($\kappa=3/2$ in $d\geq1$~\cite{LapWei-00}) For all $d \ge 1$ with $\kappa = \frac32$, we have $L_{3/2,d}=L_{3/2,d}^{\rm sc}$;
    \item ($\kappa<3/2$ is not semi-classical in 1D~\cite{LieThi-76}) For $d=1$ and $\kappa<3/2$, we have $L_{\kappa, 1} > L_{\kappa, 1}^{\rm sc}$;
    \item ($\kappa<1$ is not semi-classical~\cite{HelRob-90}) For all $d \ge 1$ and $\kappa < 1$, we have $L_{\kappa, d} > L_{\kappa, d}^{\rm sc}$;
    \item ($\kappa=0$ in $d\geq7$~\cite{GlaGroMar-78}, see also \cite{Frank-20_ppt}) We have $L_{0,d}> L_{0,d}^{\rm sc}>L_{0,d}^{(1)}$ in dimensions $d\geq 8$ and $L_{0,d}> L^{(1)}_{0,d}> L_{0,d}^{\rm sc}$ in dimension $d=7$.    
\end{itemize}
These properties imply that there is a critical number $1 \le \kappa_{\rm sc}(d) \le \frac32$ such that
\[
L_{\kappa,d} \begin{cases}
 =L_{\kappa,d}^{\rm sc}&\text{for $\kappa\geq\kappa_{\rm sc}(d)$,}\\
>L_{\kappa,d}^{\rm sc}&\text{for $\kappa<\kappa_{\rm sc}(d)$.}
             \end{cases}
\]
The exact value of $\kappa_{\rm sc}(d)$ is unknown and of course it also remains to determine what is happening below this value. 

Next we discuss the one-bound state constant $L^{(1)}_{\kappa,d}$. In Section~\ref{sec:proof_crossing} we will prove the following result, which provides some new properties of the function $\kappa\mapsto L^{(1)}_{\kappa,d}/L^{\rm sc}_{\kappa,d}$. 

\begin{theorem}[Comparing $L^{(1)}_{\kappa,d}$ with $L^{\rm sc}_{\kappa,d}$]\label{thm:crossing}
$(i)$ For every $d\geq1$, the function $\kappa\mapsto L^{(1)}_{\kappa,d}/L^{\rm sc}_{\kappa,d}$ is strictly decreasing on its interval of definition~\eqref{eq:constraint_kappa}.

\smallskip

\noindent $(ii)$ In dimensions $1\leq d\leq 7$ there is a unique $0<\kappa_{1\cap\rm sc}(d)<\infty$ such that
	$$
	\begin{cases}
		L^{(1)}_{\kappa,d} > L^{\rm sc}_{\kappa,d} & \text{if}\ \kappa<\kappa_{1\cap\rm sc}(d) \,,\\
		L^{(1)}_{\kappa,d} = L^{\rm sc}_{\kappa,d} & \text{if}\ \kappa=\kappa_{1\cap\rm sc}(d) \,,\\
		L^{(1)}_{\kappa,d} < L^{\rm sc}_{\kappa,d} & \text{if}\ \kappa>\kappa_{1\cap\rm sc}(d) \,.
	\end{cases}	
	$$

\smallskip

\noindent $(iii)$ In dimensions $d\geq 8$, one has $L^{(1)}_{\kappa,d} < L^{\rm sc}_{\kappa,d}$ for all $\kappa\geq 0$. 

\smallskip

\noindent $(iv)$ Finally, we have for $d\geq 2$,
	\begin{equation}
	 \frac{L^{(1)}_{\kappa,d}}{L^{\rm sc}_{\kappa,d}}<\frac{L^{(1)}_{\kappa,d-1}}{L^{\rm sc}_{\kappa,d-1}}
	 \qquad\text{for all}\ \kappa\geq\max\left\{0,2-\frac{d}2\right\} \,.
	 \label{eq:compare_dimensions}
	\end{equation}
  	In particular, $\kappa_{1\cap\rm sc}(d)$ is decreasing with the dimension. 
\end{theorem}

\begin{figure}[h]
\centering
\includegraphics[width=11cm]{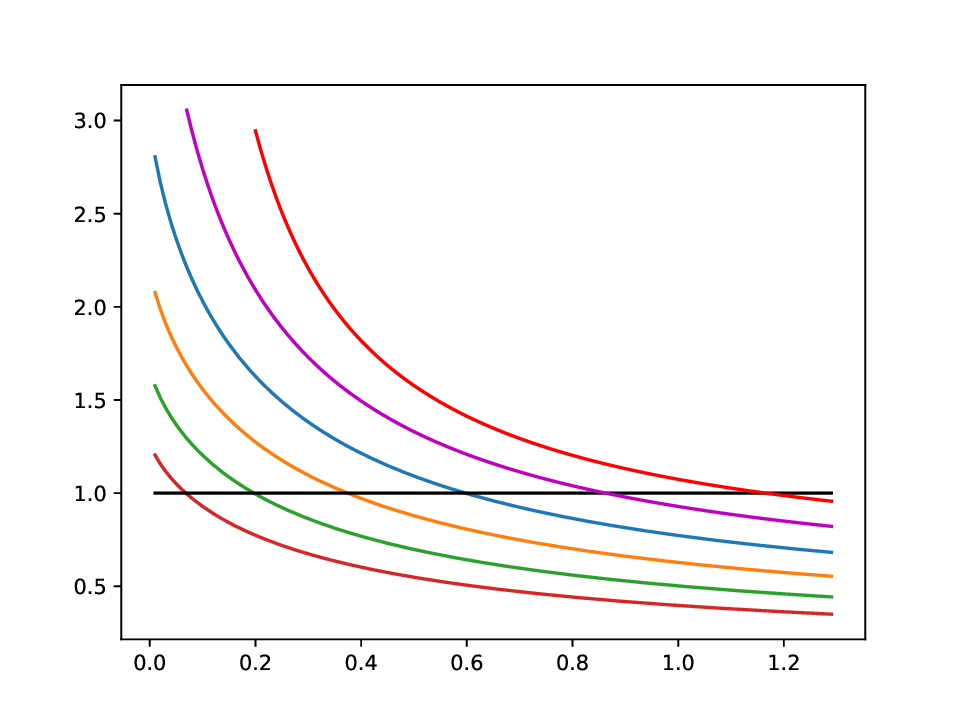}

\medskip

\begin{tabular}{|c|cccccccc|}
\hline
$d$&1&2&3&4&5&6&7&$d\geq8$\\
\hline
$\kappa_{1\cap\rm sc}(d)$&$=\frac32$&$1.1654$&$0.8627$&$0.5973$&$0.3740$&$0.1970$&$0.0683$&no crossing\\
\hline
\end{tabular}

\caption{Numerical computation of the curves $\kappa\mapsto L^{(1)}_{\kappa,d}/L^{\rm sc}_{\kappa,d}$ for $d\in\{2,...,7\}$. The curves are ordered according to the dimension, with the $d=2$ curve above the others. The points $\kappa_{1\cap\rm sc}(d)$ at which they take the value 1 are provided in the table.\label{fig:L1_over_Lsc}}
\end{figure}

That the two curves $\kappa\mapsto (L_{\kappa,d}^{(1)},L_{\kappa,d}^{\rm sc})$ cross at a unique point was part of the Lieb-Thirring conjecture \cite{LieThi-76}. In Figure~\ref{fig:L1_over_Lsc} we display a numerical computation of the curves $\kappa\mapsto L^{(1)}_{\kappa,d}/L^{\rm sc}_{\kappa,d}$ for $d\in\{2,...,7\}$ and of the crossing points $\kappa_{1\cap\rm sc}(d)$, which confirm the results of Theorem~\ref{thm:crossing}. In fact, the monotonicity with respect to the dimension~\eqref{eq:compare_dimensions} seems to hold in the whole domain of definition for $d\in\{2,3\}$. These computations complement those of Barnes in~\cite[App.~A]{LieThi-76} who only considered dimensions $d\in\{1,2,3\}$. 

The Lieb-Thirring conjecture~\eqref{eq:Lieb-Thirring_conjecture} meant that $\kappa_{\rm sc}(d)=\kappa_{1\cap\rm sc}(d)$ and that $L_{\kappa,d}=L^{(1)}_{\kappa,d}$ for $\kappa\leq \kappa_{\rm sc}(d)$. This is still believed to hold in dimension $d=1$, but not in dimensions $d\geq2$. In particular, Theorem~\ref{thm:LT_bis} implies already that 
$$\kappa_{1\cap\rm sc}(d) < \kappa_{\rm sc}(d) \qquad\text{in dimensions $2\leq d\leq 7$.}$$
The inequality is strict because otherwise we would have $L_{\kappa,d}=L_{\kappa,d}^{\rm sc}=L^{(1)}_{\kappa,d}$ at $\kappa=\kappa_{1\cap\rm sc}(d)$ which cannot hold by Theorems~\ref{thm:LT} and~\ref{thm:LT_bis}. We now discuss some further consequences of our results, mostly in the physical dimensions $d\leq3$.

\medskip

\noindent $\bullet$ \emph{In dimension $d = 1$}, we have $\kappa_{\rm sc}(1)=\kappa_{1\cap\rm sc}(1) = 3/2$. In addition, at $\kappa=1/2$, the constant is $L_{1/2,1}=L_{1/2,1}^{(1)}=1/2$ as proved in~\cite{HunLieTho-98}, with the optimal $V$ being a delta function. The remaining part of the Lieb-Thirring conjecture, namely, the equality $L_{\kappa,1}=L_{\kappa,1}^{(1)}$ for all $1/2<\kappa<3/2$, has been confirmed by numerical experiments in~\cite{Levitt-14} but it is still open.

\medskip

\noindent $\bullet$ \emph{In dimension $d = 2$}, we have $1.165\simeq\kappa_{1\cap\rm sc}(2)< \kappa_{\rm sc}(2)\leq3/2$ and this is the best we can say at present. Numerical simulations in~\cite{Levitt-14} did not provide any hint of what is happening in the region $1\leq\kappa\lesssim 1.165$. However, our Theorem~\ref{thm:LT} in dimension $d = 2$ shows that $L_{\kappa, 2} > L_{\kappa,2}^{(N)}$ for all $\kappa > 1$ and $N\geq1$. In particular, for $1 < \kappa \lesssim 1.165$, \textbf{we disprove the Lieb-Thirring conjecture that the constant is given by the $N = 1$ optimiser in 2D}. It can indeed not be given by any finite rank optimiser. 

\medskip

\noindent $\bullet$ \emph{In dimension $d = 3$}, a system with 5 bound states was numerically found in~\cite{Levitt-14} to be better than the one bound state for $\kappa\gtrsim 0.855$, showing that the one bound state case ceases to be optimal before the critical value $0.8627$ in Figure~\ref{fig:L1_over_Lsc}. Our Theorem~\ref{thm:LT_bis} implies that the one-bound state constant $L^{(1)}_{\kappa,d}$ can indeed not be optimal for all $\kappa>1/2$. This \textbf{disproves the Lieb-Thirring conjecture that the constant is given by the $N = 1$ optimiser for $1/2<\kappa\lesssim 0.8627$ in 3D.}

\medskip

\noindent $\bullet$ \emph{In dimension $d \ge 3$}, a common belief is that $\kappa_{\rm sc}(d)=1$ for all $d\geq3$. The validity of this conjecture would have some interesting physical consequences, for instance an exact lower bound involving the Thomas-Fermi kinetic energy in Density Functional Theory~\cite{LewLieSei-19_ppt}. Our Theorem~\ref{thm:LT} does not contradict this belief, since we prove that the optimal Lieb-Thirring potential cannot have a finite number of bound states. But many other situations are still possible, as we now discuss.

\medskip

Theorem~\ref{thm:LT} suggests to interpret the Lieb-Thirring inequality within the framework of statistical mechanics. For an optimal potential $V_N$ for $L^{(N)}_{\kappa,d}$, we can think of the corresponding $N$ first orthonormal eigenfunctions of $-\Delta+V_N$ as describing $N$ fermions in $\R^d$~\cite[Rmk.~8]{GonLewNaz-20_ppt}. Theorem~\ref{thm:LT} says that in the limit $N\to\ii$, the $N$ particles always attract each other, at least along a subsequence $N_k$. We \textbf{conjecture} that for $\kappa>\max\{2-d/2,0\}$ they will form a large cluster of size proportional to $N^{1/d}$ (if $\int_{\R^d} (V_N)_-^{\kappa+d/2}$ is, for instance, normalised to $N$) and that $V_N$ will converge in the limit to a bounded, but non-decaying potential $V_\ii$. There would then be no optimiser for the Lieb-Thirring constant $L_{\kappa,d}$. The semi-classical constant $L_{\kappa,d}^{\rm sc}$ corresponds to the case where the limiting potential $V_\ii$ is constant over $\R^d$, that is, the system is translation-invariant. In statistical mechanics, this is called a \emph{fluid phase}. In principle, the limiting potential $V_\ii$ could also be a non-trivial periodic function, which is then interpreted as a \emph{solid phase}. We see no obvious physical reasons for discarding this possibility, in particular in low dimensions where periodic systems are ubiquitous~\cite{BlaLew-15}. This mechanism does not seem to have been considered before in the context of Lieb-Thirring inequalities. In particular, it seems natural to conjecture that the system is in a solid phase for all $2-d/2<\kappa<\kappa_{\rm sc}(d)$ in dimensions $d=2,3$. In~\cite{FraGonLew-20b_ppt} we shall discuss this new point of view in detail. 

\begin{remark}
In dimension $d = 2$, some preliminary numerical tests suggest that the difference $L_{\kappa, 2} - L_{\kappa, 2}^{(1)}$ might be very small in the region $1 < \kappa \lesssim 1.165$. This makes the problem difficult to simulate as we need high precision. 
\end{remark}

\subsection{Dual Lieb-Thirring inequalities}\label{sec:duallt}
Our strategy to prove Theorem~\ref{thm:LT} is to study the dual version of the Lieb-Thirring inequality~\eqref{eq:LT_V_N}. This dual version is well known for $\kappa = 1$ and it is often used in practical applications. The dual inequality for $\kappa>1$ appears, for instance, in~\cite{LioPau-93}, but is less known and we briefly recall it in this subsection. There is no known dual problem for $\kappa<1$, except for a certain substitute for $\kappa=0$ in dimensions $d\geq3$~\cite{Frank-14}.

Let $0 \le \gamma  = \gamma^*$ be a self-adjoint non-negative operator of $\rank(\gamma)\leq N$, of the form $\gamma = \sum_{j=1}^N n_j | u_j \ket \bra u_j |$ with $u_1,...,u_N$ an orthonormal family in $L^2(\R^d)$. For $1 \le q < \infty$, we denote by
\[
    \|\gamma\|_{\gS^q}:=(\tr|\gamma|^q)^{1/q} = \left( \sum_{j=1}^N n_j^q \right)^{1/q}
\]
its $q$-th Schatten norm~\cite{Simon-79}, and use the convention that $\| \gamma \|_{\gS^\ii} = \| \gamma \|$ is the operator norm. The density of $\gamma$ is the function $\rho_\gamma \in L^1(\R^d)$ defined by
\[
    \rho_\gamma(x) := \sum_{j=1}^N n_j | u_j (x)|^2,
\]
and the kinetic energy of $\gamma$ is 
\[
    \Tr( - \Delta \gamma) := \sum_{j=1}^N n_j \int_{\R^d} | \nabla u_j |^2 (x) \rd x
\]
with the convention that $\Tr( - \Delta \gamma) = + \infty$ if $u_j \notin H^1(\R^d)$ for some $j$. Let $1 \le p \le 1 + \frac{2}{d}$ with $d \ge 1$, and let
\begin{equation*}
q := \begin{cases}
\frac{2p+d-dp}{2+d-dp}&\text{for $1\le p<1+\frac2d$,}\\
+\ii&\text{for $p=1+\frac2d$.}
\end{cases}
\end{equation*}
We denote by $K_{p, d}^{(N)}$ the best (that is, largest possible) constant in the inequality
\begin{equation}
 \boxed{ K_{p,d}^{(N)} \norm{\rho_\gamma}_{L^p(\R^d)}^{\frac{2p}{d(p-1)}} \le  \norm{\gamma}_{\gS^q}^{\frac{p(2-d)+d}{d(p-1)}}\tr(-\Delta\gamma)}
 \label{eq:LT_Schatten}
\end{equation}
valid for all $0 \le \gamma = \gamma^*$ with $\rank(\gamma) \le N$. The fact that $K_{p, d}^{(N)}$ is well-defined with $K_{p, d}^{(N)} > 0$ is a consequence of the next result, together with the Lieb-Thirring theorem.

\begin{lemma}[Duality]\label{lem:duality_N}
    Let $1\leq N\leq \ii$, $d\geq1$ and $1\leq p\leq 1+\frac2d$, and set
    \[
        \kappa := \frac{p}{p-1} - \frac{d}{2}, \quad \text{so that} \quad \frac{\kappa}{\kappa - 1} = q.
    \]
    Then,
    \begin{equation} \label{eq:duality_K_and_L}
    K_{p,d}^{(N)} \left( L_{\kappa,d}^{(N)}\right)^{\frac{2}d}=\left(\frac{\kappa}{\kappa+\frac{d}{2}}\right)^{\frac{2\kappa}d}\left(\frac{d}{2\kappa + d}\right).
    \end{equation}
\end{lemma}

The lemma says that the inequality~\eqref{eq:LT_Schatten} is dual to the finite-rank Lieb-Thirring inequality~\eqref{eq:LT_V_N}. This is because the density $\rho_\gamma$ is the variable dual to the potential $V$ whereas the density matrix $\gamma$ can be interpreted as the dual of the Schr\"odinger operator $-\Delta+V$. Hence $p$ is the dual exponent of $\kappa+d/2$ and $q$ the one of $\kappa$. 
The proof of Lemma~\ref{lem:duality_N}, provided in Appendix~\ref{appendix:proof_duality}, also shows how to relate the corresponding optimisers, assuming they exist. A similar argument, but without the constraint on the rank, can be found for instance in~\cite{LioPau-93}.

We denote
$$
K_{p,d} := \lim_{N\to\infty} K_{p,d}^{(N)} = \inf_{N\geq1} K_{p,d}^{(N)} \,.
$$
This constant is related to the constant $L_{\kappa,d}$ in \eqref{eq:LT_V_form} by
\begin{equation} \label{eq:duality_K_and_L_unconst}
K_{p,d} \left( L_{\kappa,d} \right)^{\frac{2}d}=\left(\frac{\kappa}{\kappa+\frac{d}{2}}\right)^{\frac{2\kappa}d}\left(\frac{d}{2\kappa + d}\right)
\end{equation}
and is the best constant in the inequality
\begin{equation}
\boxed{ K_{p,d} \norm{\rho_\gamma}_{L^p(\R^d)}^{\frac{2p}{d(p-1)}} \le  \norm{\gamma}_{\gS^q}^{\frac{p(2-d)+d}{d(p-1)}}\tr(-\Delta\gamma)}
\label{eq:LT_Schatten_unconst}
\end{equation}
valid for all $0 \le \gamma = \gamma^*$.

In Section~\ref{sec:proof_K}, we study the dual problem \eqref{eq:LT_Schatten} and prove the following result which, together with Lemma~\ref{lem:duality_N}, immediately implies Theorem~\ref{thm:LT}.

\begin{theorem}[Existence of optimisers and properties]
    \label{thm:K}
    Let $d\geq1$ and $1\leq p\leq 1+2/d$.
    
    \medskip
    
    \noindent $(i)$ \textbf{Existence.} For every finite $N\geq1$, the problem $K^{(N)}_{p,d}$ in~\eqref{eq:LT_Schatten} admits an optimiser $\gamma$. 
    
    \medskip
    
    \noindent $(ii)$ \textbf{Equation.} After an appropriate normalisation, any optimiser $\gamma$ for $K_{p,d}^{(N)}$ has rank $1\leq R\leq N<\ii$ and can be written in the form 
    $$\gamma=\sum_{j=1}^Rn_j|u_j\rangle\langle u_j|$$
with
\begin{equation}
 n_j=\begin{cases}
    \left(\frac{2p}{d(p-1)}\right)^{\frac{1}{p-1}}\frac{2p+d-dp}{d(p-1)}\frac{|\mu_j|^{\frac{1}{q-1}}}{\sum_{k=1}^R|\mu_k|^{\frac{q}{q-1}}}&\text{for $p<1+\frac2d$,}\\
    \frac2d \left(\frac{d}{d+2}\right)^{\frac{1}{p-1}}\frac{1}{\sum_{k=1}^R|\mu_k|}&\text{for $p=1+\frac2d$,}
    \end{cases}
    \label{eq:formula_n_j}
\end{equation}
    where the corresponding orthonormal system $(u_1,...,u_R)$ solves the nonlinear Schr\"o\-din\-ger equation
    \begin{equation} \label{eq:NLS_in_Lemma}
    \forall j = 1, \cdots, R, \quad \Big(-\Delta-\rho_\gamma(x)^{p-1}\Big)u_j=\mu_j\,u_j,
    \quad \text{with} \quad
    \rho_\gamma = \sum_{j=1}^R n_j | u_j |^2.
    \end{equation}
    Here $\mu_j$ are the $R$ first negative eigenvalues of $H_\gamma := - \Delta - \rho_\gamma^{p-1}$. In particular, this operator has at least $R$ negative eigenvalues. If $R<N$, then it has exactly $R$ negative eigenvalues. Finally, the potential $V=-\rho_\gamma^{p-1}$ is an optimiser for the finite-rank Lieb-Thirring problem $L^{(N)}_{\kappa,d}$ in~\eqref{eq:LT_V_N}.
    
    \medskip
    
    \noindent $(iii)$ \textbf{Rank.} If, in addition, $p<2$, then there exists an infinite sequence of integers $N_1=1<N_2=2<N_3<\cdots$ so that 
    $$K^{(N_{k})}_{p,d} < K^{(N_k-1)}_{p,d}$$
    and any optimiser for $K^{(N_k)}_{p,d}$ must have rank $R=N_k$. In particular,
    $$K_{p,d}<K_{p,d}^{(N)},\qquad \text{for all} \ N\geq1.$$
\end{theorem}

The assertions in $(i)$ and $(ii)$ follow by applying well-known methods from the calculus of variation adapted to the setting of operators; see, for instance, \cite{Solovej-91,Bach-93,FraLieSeiSie-07,Lewin-11}. 
For $(iii)$, we use ideas from~\cite{GonLewNaz-20_ppt}, which consist in evaluating the exponentially small interaction between two copies of an optimiser placed far from each other, in order to show that
$$K^{(2N)}_{p,d}<K^{(N)}_{p,d}$$
whenever $K^{(N)}_{p,d}$ admits an optimiser of rank $N$. The proof is provided in Section~\ref{sec:proof_K} below. This argument inspired our proof of Theorem~\ref{thm:LT_bis} for $\kappa<1$ and $N=2$, provided in Section~\ref{sec:proof_LT_V_kappa<1}. There we use the $N=1$ Gagliardo-Nirenberg optimiser to construct a trial state for $N=2$ but we do not prove the existence of an optimal potential.

\subsection{Fermionic Nonlinear Schr\"odinger Equation}
The system of coupled nonlinear equations~\eqref{eq:NLS_in_Lemma} has some similarities with that studied in~\cite{GonLewNaz-20_ppt}, where one has $n_j=1$ instead of~\eqref{eq:formula_n_j}. Here we exhibit a link between the two problems and use this to solve a question left open in~\cite{GonLewNaz-20_ppt}.

In~\cite{GonLewNaz-20_ppt} the authors studied the minimisation problem 
\begin{equation}
\boxed{J(N)=\inf \left\{ \tr(-\Delta\gamma)-\frac1p\int_{\R^d}\rho_\gamma(x)^p\,\rd x : \ 0\leq \gamma=\gamma^*\leq1,\ \Tr(\gamma) = N \right\}.}
\label{eq:def_J1}
\end{equation}
Under the assumption $1 < p < 1 + {2}/{d}$, it is proved in~\cite{GonLewNaz-20_ppt} that $- \infty < J(N) < 0$ for all $N > 0$. Under the additional assumption that $p < 2$, it was also shown that there is an infinite sequence of integers $N_1 = 1 < N_2 = 2 < N_3 < \cdots$ such that $J(N_k)$ has a minimiser $\gamma$ of rank $N_k$. This minimiser is a projector of the form $\gamma = \sum_{j=1}^{N_k} | u_j \ket \bra u_j|$, where $u_1,...,u_{N_k}$ form an orthonormal system and solve the \emph{fermionic NLS equation}
\begin{equation}
\forall j = 1, \cdots, N_k, \qquad \left(-\Delta- \rho_\gamma(x)^{p-1}\right)u_j=\mu_j\,u_j, \quad \text{with} \quad 
\rho_\gamma = \sum_{i=1}^{N_k}|u_i|^2.
\label{eq:fermionic_NLS}
\end{equation}
Here again $\mu_1 < \mu_2 \le \cdots \le \mu_{N_k} < 0$ are the $N_k$ first eigenvalues of $H_\gamma := - \Delta - \rho_\gamma^{(p-1)}$. The existence of minimisers for $J(N_k)$ therefore proves the existence of solutions of the fermionic NLS equation~\eqref{eq:fermionic_NLS}, for all $1 \le p < \min\{ 2, 1 + 2/d\}$ and $N = N_k$. In dimension $d = 1$, this does not cover the case $p \in [2, 3)$. In the present paper, we prove the following result for the case $p = 2$, which was announced in~\cite{GonLewNaz-20_ppt} and actually also follows from the analysis in~\cite{LieLla-78}.

\begin{theorem}[Non-existence of minimisers for $d=1$, $p=2$]\label{th:nonExistence_J(N)}
    Let $d=1$ and $p=2$. For all $N \ge 1$, we have $J(N)=N\, J(1)$. In addition, for all $N \ge 2$, $J(N)$ admits no minimiser.
\end{theorem}

The theorem is reminiscent of a similar result for the true Schr\"odinger (Lieb-Liniger~\cite{LieLin-63}) model in 1D describing $N$ particles interacting with the delta potential. In the attractive case, only two-particle (singlet) bound states exist~\cite{McGuire-64,Yang-68,LieLla-78}. The same result in the Hartree-Fock case was proved in~\cite{LieLla-78}. The spatial component of the singlet state coincides with our $N=1$ solution. 

In the case $N = 1$ and $1<p<1+2/d$, it is proved in~\cite[Lem.~11]{GonLewNaz-20_ppt} that $J(1)$ has the Gagliardo-Nirenberg-Sobolev optimiser $\gamma = | U \ket \bra U |$, where 
\begin{equation}
 U(x)=m^{-\frac{p-1}{2(1+2/d-p)}-\frac12}\;Q\left(m^{-\frac{p-1}{d(1+2/d-p)}}x\right),\qquad\int_{\R^d}U(x)^2\,{\rm d}x=1,
 \label{eq:sol_NLS_mass_1}
\end{equation}
and $Q$ is the unique positive radial solution to the NLS equation
\begin{equation}
-\Delta Q-Q^{2p-1}+Q =0, \quad \text{with mass} \quad m := \int_{\R} Q^2.
\label{eq:NLS}
\end{equation}
When $d=1$ and $p=2$, we have the explicit formula
$$U(x)=\frac{1}{2^{\frac32}\cosh(x/4)}.$$
Our strategy to prove Theorem~\ref{th:nonExistence_J(N)} for $d=1$ is to relate $J(N)$ to the dual Lieb-Thirring constant $K_{\kappa, 1}^{(N)}$ for $\kappa=3/2$, and use $K_{3/2, 1}^{(N)} = K_{3/2, 1}^{(1)}$. The proof is given in Section~\ref{ssec:proof_nonExistence_J(N)} below.

The same argument gives that if the Lieb-Thirring conjecture $K_{\kappa, 1}^{(N)} = K_{\kappa, 1}^{(1)}$ is true for some $1<\kappa<3/2$, then $J(N)=N\,J(1)$ for $p=(\kappa+1/2)/(\kappa-1/2)$; see Remark \ref{rem:ltconjj}.

\medskip

Even if $J(N)$ has no minimiser for $N\geq 2$ if $d=1$ and $p=2$, one may still wonder whether the fermionic NLS equation~\eqref{eq:fermionic_NLS} possesses orthonormal solutions. We believe there are no other solutions than the $N=1$ case and are able to prove this for $N=2$, using the fundamental fact that the system is completely integrable~\cite{Manakov-74}. The following is stronger than Theorem~\ref{th:nonExistence_J(N)} for $N=2$.

\begin{theorem}[Non-existence of solutions for $p = 2$, $d=1$ and $N = 2$] \label{th:N=2}
    Let $\mu_1 \le \mu_2 < 0$, and let $u_1, u_2$ be two square integrable real-valued functions solving 
    \begin{equation} \label{eq:no_binding_u}
    \begin{cases}
    - u_1'' - (u_1^2 + u_2^2) u_1 = \mu_1 u_1, \\
    - u_2'' - (u_1^2 + u_2^2) u_2 = \mu_2 u_2.
    \end{cases}
    \end{equation}
    If $\| u_1 \|_{L^2(\R)} = \| u_2 \|_{L^2(\R)}=1$, then we have $\mu_1 = \mu_2$ and 
    \begin{equation}
     u_1(x)=\pm \frac{1}{2\cosh\big((x-x_0)/2\big)},\qquad u_2(x)=\pm \frac{1}{2\cosh\big((x-x_0)/2\big)}
     \label{eq:formulas_u_1_u_2}
    \end{equation}
    for some $x_0\in\R$ and two uncorrelated signs $\pm$.
\end{theorem}

The proof can probably be generalised to show that there are no solutions for all $N\geq3$ at $p=2$ but we only address the simpler case $N=2$ here. The proof is given in Section~\ref{ssec:proof_N=2}. More comments about the NLS problem~\eqref{eq:def_J1} can be read in Appendix~\ref{app:NLS_comments}.

\subsection*{Structure of the paper}
In the next section we recall useful facts about the Lieb-Thirring constant $L^{(1)}_{\kappa,d}$ and provide the proof of Theorem~\ref{thm:crossing}. In 
Section~\ref{sec:proof_K}, we prove Theorem~\ref{thm:K}, which implies Theorem~\ref{thm:LT}. Section~\ref{sec:proof_LT_V_kappa<1} is devoted to the proof of Theorem~\ref{thm:LT_bis}. We prove Theorem~\ref{th:nonExistence_J(N)} and Theorem~\ref{th:N=2} in Sections~\ref{ssec:proof_nonExistence_J(N)} and~\ref{ssec:proof_N=2}, respectively. The proof of duality (Lemma~\ref{lem:duality_N}) is given in Appendix~\ref{appendix:proof_duality} whereas Appendix~\ref{app:NLS_comments} contains more comments on the NLS model from~\cite{GonLewNaz-20_ppt}. Finally, in Appendix~\ref{app:LT-Sobolev} we compare our results with those in~\cite{HonKwoYoo-19}.


\section{The one-bound state constant $L^{(1)}_{\kappa,d}$: Proof of Theorem~\ref{thm:crossing}}
\label{sec:proof_crossing}

In this section we discuss some properties of the one-bound state constant $L^{(1)}_{\kappa,d}$ and provide the proof of Theorem~\ref{thm:crossing}. The Gagliardo-Nirenberg inequality states that
\begin{equation}
K_{p,d}^{\rm GN}\left(\int_{\R^d}|u(x)|^{2p}\,\rd x\right)^{\frac2{d(p-1)}}\leq \left(\int_{\R^d}|\nabla u(x)|^2\,\rd x\right)\left(\int_{\R^d}|u(x)|^2\,\rd x\right)^{\frac{(2-d)p+d}{d(p-1)}}
\label{eq:GN}
\end{equation}
for all
\begin{equation*}
\begin{cases}
1< p <+\ii&\text{for $d=1,2$,}\\
1 < p \leq\frac{d}{d-2}&\text{for $d\geq3$,}
\end{cases}
\end{equation*}
with the best constant $K_{p,d}^{\rm GN}>0$. In dimension $d=1$ one can take $p\to+\ii$. The constants $K_{p,1}^{\rm GN}$ and the optimisers are known explicitly in $d=1$ \cite{Nagy-41}. In particular, the optimiser is unique up to translations, dilations and multiplication by a phase factor. As explained, for instance, in~\cite{Tao-06,Frank-13,CarFraLie-14}, by combining the results on existence~\cite{Strauss-77,BerLio-83,Weinstein-83}, symmetry~\cite{GidNiNir-81,AlvLioTro-86} and uniqueness~\cite{Coffman-72,Kwong-89,McLeod-93} one infers that in any $d\geq 2$ as well, there is a unique optimiser $Q$, up to translations, dilations and multiplication by a phase factor. This function can be chosen positive and to satisfy~\eqref{eq:NLS} when $p<1+2/d$. When $p=1+2/d$, it still can be chosen positive and to satisfy the equation in~\eqref{eq:NLS}. The integral $\int_{\R^d} Q^2\,\rd x$ will be a dimension-dependent constant.

For an operator $\gamma$ of rank one the inequality~\eqref{eq:LT_Schatten} is equivalent to~\eqref{eq:GN}, hence we obtain
\begin{equation}
 K_{p,d}^{(1)}=K_{p,d}^{\rm GN}.
\label{eq:K_1}
\end{equation}
The duality argument from Lemma~\ref{lem:duality_N} shows that 
\begin{equation}
L_{\kappa,d}^{(1)}=\left(\frac{2\kappa}{2\kappa+d}\right)^{\kappa+\frac{d}2}\left(\frac{d}{2\kappa}\right)^{\frac{d}2}\left(K_{p,d}^{\rm GN}\right)^{-\frac{d}2}<\ii.
\label{eq:formula_L_1}
\end{equation}
By the implicit function theorem and the non-degeneracy of $Q$~\cite{McLeod-93,Tao-06,Frank-13}, the Gagliardo-Nirenberg constant $K_{p,d}^{(1)}$ is known to be real-analytic in $p$, so that $L_{\kappa,d}^{(1)}$ is a real-analytic function of $\kappa$. In this paper we will only use the continuity of $\kappa\mapsto L_{\kappa,d}^{(1)}$, which is more elementary and which we explain now for completeness. We claim that $p\mapsto K_{p,d}^{(1)}$ is continuous in the interval $(1,\infty)$ if $d=1,2$ and $(1,d/(d-2)]$ if $d\geq 3$, which implies the continuity of $L^{(1)}_{\kappa,d}$ on the corresponding intervals. To prove this fact, we can notice that 
\begin{multline}
\log\left(\big(K_{p,d}^{(1)}\big)^{\frac{d(p-1)}{4p}}\right)=\inf_{u\in H^1(\R^d)}\bigg\{\frac{d(1-p^{-1})}{2}\log\norm{\nabla u}_{L^2(\R^d)}\\
+\left(1-\frac{d(1-p^{-1})}{2}\right)\log\norm{u}_{L^2(\R^d)}-\log\norm{u}_{L^{2p}(\R^d)}\bigg\}.
\end{multline}
By H\"older's inequality, $p^{-1}\mapsto \log\norm{u}_{L^{2p}(\R^d)}$ is convex. Hence after minimising over $u$ we find that $p\mapsto \big(K_{p,d}^{(1)}\big)^{\frac{d(p-1)}{4p}}$ is upper semi-continuous on $[1,\infty)$ if $d=1,2$ and on $[1,d/(d-2)]$ if $d\geq 3$, and log-concave in $1/p$. Log-concavity implies continuity on the interior of the interval of definition and then upper semicontinuity implies continuity up to the endpoints.

Our goal in the rest of this section is to compare $L^{(1)}_{\kappa,d}$ with the semi-classical constant $L^{\rm sc}_{\kappa,d}$.  First, the argument from~\cite{AizLie-78} can be used to prove that $\kappa\mapsto L^{(1)}_{\kappa,d}/L^{\rm sc}_{\kappa,d}$ is non-increasing. We show here that it is even strictly decreasing, which is $(i)$ in Theorem~\ref{thm:crossing}. 

\begin{lemma}\label{al}
	For any $d\geq 1$, the function $\kappa\mapsto L^{(1)}_{\kappa,d}/L_{\kappa,d}^{\rm sc}$ is strictly decreasing.
\end{lemma}

\begin{proof}
	Following~\cite{AizLie-78}, we use the fact that for any $0\leq\kappa'<\kappa$ and $\lambda\in\R$, we have
	\begin{eqnarray}
	\label{eq:betafcn}
	\lambda_-^\kappa = c_{\kappa,\kappa'} \int_0^\infty (\lambda+t)_-^{\kappa'} \,t^{\kappa-\kappa'-1}\,\rd t
	\end{eqnarray}
	for some constant $c_{\kappa,\kappa'}>0$. Let $V\in L^{\kappa+d/2}(\R^d)$. By the variational principle we have $(\lambda_1(-\Delta+V)+ t)_- \leq  |\lambda_1(-\Delta-(V+t)_-) |$ for any $t\geq 0$
	and we can bound, using the definition of $L_{\kappa',d}^{(1)}$,
	\begin{align}
	\label{eq:alproof}
	\big(\lambda_1(-\Delta+V)+ t\big)_-^{\kappa'} &\leq \big|\lambda_1\big(-\Delta-(V+t)_-\big)\big|^{\kappa'}\notag\\
	& \leq L_{\kappa',d}^{(1)} \int_{\R^d} \left( V(x)+t \right)_-^{\kappa'+\frac d2} \rd x \notag \\
	& = L_{\kappa',d}^{(1)} \left( L_{\kappa',d}^{\rm sc} \right)^{-1} \iint_{\R^d\times\R^d} \left( |\xi|^2 + V(x)+t \right)_-^{\kappa'} \frac{\rd\xi\,\rd x}{(2\pi)^d} \,.
	\end{align}
	Thus, integrating over $t$ using \eqref{eq:betafcn} on both sides, we obtain
	\begin{align}
	\label{eq:alproof1}
	\lambda_1(-\Delta+V)_-^\kappa & \leq L_{\kappa',d}^{(1)} \left( L_{\kappa',d}^{\rm sc} \right)^{-1} \iint_{\R^d\times\R^d} \left( |\xi|^2 + V(x) \right)_-^{\kappa} \frac{\rd\xi\,\rd x}{(2\pi)^d} \notag \\
	& = L_{\kappa',d}^{(1)} \left( L_{\kappa',d}^{\rm sc} \right)^{-1} L_{\kappa,d}^{\rm sc} \int_{\R^d} V(x)_-^{\kappa+\frac d2}\,\rd x \,.
	\end{align}
	This shows that
	\begin{equation}
	\label{eq:alproof2}
	L_{\kappa,d}^{(1)} \leq L_{\kappa',d}^{(1)} \left( L_{\kappa',d}^{\rm sc} \right)^{-1} L_{\kappa,d}^{\rm sc} \,,
	\end{equation}
	that is, $\kappa\mapsto L^{(1)}_{\kappa,d}/L_{\kappa,d}^{\rm sc}$ is nonincreasing.	
	
	As was recalled at the beginning of this section, it is known that for the optimisation problem corresponding to $L_{\kappa',d}^{(1)}$ there is an optimiser. This optimiser is a power of the solution of the positive solution of \eqref{eq:NLS} and therefore does not vanish. Since for any $V\in L^{\kappa+d/2}(\R^d)$ and for any $t>0$, the function $-(V+t)_-$ is supported on a set of finite measure, this function cannot be an optimiser for $L_{\kappa',d}^{(1)}$. Therefore the second inequality in \eqref{eq:alproof} is strict for all $t>0$ and, consequently, inequality \eqref{eq:alproof1} is strict for any $V\in L^{\kappa+d/2}(\R^d)$. Taking, in particular, $V$ to be an optimiser corresponding to $L_{\kappa,d}^{(1)}$, we obtain that inequality \eqref{eq:alproof2} is strict, which is the assertion of the lemma.	
\end{proof}

Next, we prove an inequality relating the constant $L^{(1)}_{\kappa,d}$ with the ones in lower dimensions, in the spirit of the Laptev-Weidl method of lifting dimensions~\cite{LapWei-00}. 

\begin{lemma}\label{lem:lifting_dimensions}
For any $d\geq2$ and $\kappa>0$, we have
\begin{equation}
L^{(1)}_{\kappa,d}<L^{(1)}_{\kappa,d-n}\,L^{(1)}_{\kappa+\frac{d-n}{2},n},\qquad \forall n\in\{1,...,d-1\}.
\label{eq:lifting_dimensions}
\end{equation}
The same inequality holds for $\kappa=0$ if $d-n\geq 3$.
\end{lemma}

\begin{proof}
Let $V$ be the optimizer for $L^{(1)}_{\kappa,d}$ with corresponding ground state $u$, which can both be expressed in terms of the NLS solution $Q$ in~\eqref{eq:NLS}. We write $x=(x_1,x_2)\in\R^{d-n}\times\R^n$ and denote by $\lambda(x_1)$ the first eigenvalue of $-\Delta_{x_2}+V(x_1,\cdot)$ in $\R^n$. Writing 
$$-\Delta+V=-\Delta_{x_1}+\big(-\Delta_{x_2}+V(x_1,x_2)\big)\geq -\Delta_{x_1}+\lambda(x_1)$$ 
and taking the scalar product with $u$, we find $\lambda_1\left(-\Delta+V\right)> \lambda_1\left(-\Delta_{x_1}+\lambda(x_1)\right)$. The strict inequality is because $u$ does not solve an eigenvalue equation in $x_1$ at fixed $x_2$. This gives
\begin{align*}
\left|\lambda_1(-\Delta+V) \right|^\kappa&<\left|\lambda_1(-\Delta_{x_1}+\lambda(x_1)) \right|^\kappa\\
&\leq L^{(1)}_{\kappa,d-n}\int_{\R^n}\left|\lambda(x_1)\right|^{\kappa+\frac{d-n}{2}}\,\rd x_1\\
&\leq L^{(1)}_{\kappa,d-n}L^{(1)}_{\kappa+\frac{d-n}{2},n}\iint_{\R^{d-n}\times\R^n}V(x_1,x_2)_-^{\kappa+\frac{d}2}\,\rd x_1\,\rd x_2
\end{align*}
and we obtain inequality~\eqref{eq:lifting_dimensions} for $\kappa>0$.

The proof for $\kappa=0$ if $d-n\geq 3$ is similar. Let again $V$ be the optimizer for $L^{(1)}_{\kappa,d}$ and $u$ the corresponding ground state. More precisely, $u$ is an eigenfunction corresponding to the eigenvalue zero if $d\geq 5$ and it is a zero energy resonance function (that is, an element of $\dot H^1(\R^d)\setminus L^2(\R^d)$) if $d=3,4$. We have
\begin{align*}
0 & = \int_{\R^d} \left( |\nabla u|^2 + V|u|^2\right)\rd x \geq \int_{\R^d} \left( |\nabla_{x_1}u|^2 + \lambda(x_1)|u|^2 \right)\rd x \\
& \geq \int_{\R^n}  \bigg\{ \int_{\R^{d-n}} |\nabla_{x_1} u|^2\,\rd x_1\\
&\qquad \qquad- \left( \int_{\R^{d-n}} |\lambda(x_1)|^{\frac{d-n}{2}}\,\rd x_1 \right)^\frac{2}{d-n} \left( \int_{\R^{d-n}} |u|^\frac{2(d-n)}{d-n-2}\,\rd x_1 \right)^\frac{d-n-2}{d-n} \bigg\} \rd x_2 \\
& \geq \int_{\R^d} |\nabla_{x_1}u|^2\,\rd x  \left( 1 - S_{d-n}^{-1} \left( \int_{\R^{d-n}} |\lambda(x_1)|^{\frac{d-n}{2}}\,\rd x_1 \right)^\frac{2}{d-n} \right),
\end{align*}
where $S_{d-n}$ is the optimal Sobolev constant in dimension $d-n$. We conclude that
$$
\int_{\R^{d-n}} |\lambda(x_1)|^{\frac{d-n}{2}}\,\rd x_1 \geq S_{d-n}^\frac{d-n}{2} = \left( L_{0,d-n}^{(1)} \right)^{-1} .
$$
On the other hand, we have
$$
|\lambda(x_1)|^\frac{d-n}{2} < L_{\kappa+\frac{d-n}{2},n}^{(1)} \int_{\R^n} |V(x_1,x_2)|^\frac{d}{2}\,\rd x_1 \,,
$$
where the strict inequality follows from the fact that for no $x_1$, $V(x_1,\cdot)$ is an optimal potential for $L_{\kappa+\frac{d-n}{2},n}^{(1)}$. (Indeed, $V(x_1,\cdot)$ is algebraically decaying, whereas we know that the optimal potential for $L_{\kappa+\frac{d-n}{2},n}^{(1)}$ is exponentially decaying.) Combining the last two inequalities we obtain
$$
 L_{\kappa+\frac{d-n}{2},n}^{(1)} \int_{\R^d} |V(x)|^\frac{d}{2}\,\rd x > \left( L_{0,d-n}^{(1)} \right)^{-1} \,.
$$
Since $\int_{\R^d} |V(x)|^\frac{d}{2}\,\rd x=(L_{0,d}^{(1)})^{-1}$, this is the claimed inequality for $\kappa=0$.
\end{proof}

Note that the semi-classical constants satisfy the relation
\begin{equation}
L^{\rm sc}_{\kappa,d}=L^{\rm sc}_{\kappa,d-n}\,L^{\rm sc}_{\kappa+\frac{d-n}{2},n},\qquad \forall n\in\{1,...,d-1\}
\label{eq:lifting_dimensions_sc}
\end{equation}
so that we obtain the same inequality as~\eqref{eq:lifting_dimensions} for $L^{(1)}_{\kappa,d}/L^{\rm sc}_{\kappa,d}$. According to the Laptev-Weidl method of lifting dimensions the bound~\eqref{eq:lifting_dimensions} with $\leq$ instead of $<$ holds for the matrix-valued Lieb-Thirring constants $L_{\kappa,d}^{({\rm mat})}$. Using the results from~\cite{LapWei-00,HunLapWei-00}, one sees that for $n=d-1$ and $\kappa\in\{1/2\}\cup[3/2,\ii)$ the bound~\eqref{eq:lifting_dimensions} with $\leq$ instead of $<$ holds for the usual Lieb-Thirring constants $L_{\kappa,d}$. One might wonder whether this is true more generally.

In~\cite{Martin-90}, Martin used a similar idea but instead of removing he added one dimension by considering the potential $W(x,t):=V(x)+\lambda t^2$. This led to the inequality
$$\frac{L_{\kappa',d}^{(N)}}{L^{\rm sc}_{\kappa',d}}< \frac{L^{(N)}_{\kappa,d+1}}{L^{\rm sc}_{\kappa,d+1}},\qquad\forall \kappa'\geq\kappa+\frac{1}2$$
for all $1\leq N\leq \ii$. This can in fact be improved for $N=1$, see~\cite[Sec.~3]{Martin-90}. 

The proof of Theorem~\ref{thm:crossing} follows from Lemmas~\ref{al} and~\ref{lem:lifting_dimensions}. 

\begin{proof}[Proof of Theorem~\ref{thm:crossing}]
In $d=1$, the constant $L^{(1)}_{\kappa,1}$ is explicit and the unique intersection at $\kappa_{1\cap\rm sc}(1)=3/2$ follows by explicit comparison. The bound~\eqref{eq:compare_dimensions}, for general $d\geq 2$, then follows immediately from~\eqref{eq:lifting_dimensions} with $n=1$, by using~\eqref{eq:lifting_dimensions_sc} and the fact that $L^{(1)}_{\kappa,1}> L^{\rm sc}_{\kappa,1}$ for $\kappa>3/2$.

In dimension $d\geq 3$ using the explicit formula for the sharp Sobolev constant~\cite{Rodemich-66,Aubin-76,Talenti-76} (see also~\cite[Thm.~8.3]{LieLos-01}) we obtain the exact formula at $\kappa=0$:
\begin{equation}
\frac{L^{(1)}_{0,d}}{L^{\rm sc}_{0,d}}=2^{d-1}d^{-\frac{d}2}(d-2)^{-\frac{d}2} d!
\label{eq:L1_Lsc_kappa_0}
\end{equation}
This is larger than 1 in dimensions $d\in\{3,...,7\}$ but smaller than 1 in dimension $d\geq8$, as noted in~\cite{LieThi-76,GlaGroMar-78}. In fact, this is decreasing with the dimension for $d\geq4$ by~\eqref{eq:compare_dimensions} and the value in dimension $d=8$ equals $L^{(1)}_{0,8}/L^{\rm sc}_{0,8}\simeq0.9722$. Thus, if $d\geq 8$ the part $(iii)$ of the theorem follows from Lemma~\ref{al}.

In dimension $d=2$, simple numerical computations provide $L^{(1)}_{1,2}/L_{1,2}^{\rm sc}\simeq 1.074>1$ at $\kappa=1$, see~\cite{LieThi-76,Weinstein-83}. Alternatively, to see this analytically, one can use the trial function $u(x)=e^{-|x|^2}$ in the Gagliardo-Nirenberg inequality \eqref{eq:GN} to obtain an upper bound  on the constant $K_{2,2}^{\rm GN}=K_{2,2}^{(1)}$. Via \eqref{eq:formula_L_1} this gives the lower bound $L_{1,2}^{(1)}\geq (8\pi)^{-1} = L_{1,2}^{\rm sc}$. Since the Gaussian does not satisfy the Euler-Lagrange equation for $K_{2,2}^{\rm GN}$, the inequality is, in fact, strict, as claimed.

On the other hand, it also follows from~\eqref{eq:compare_dimensions} that $L^{(1)}_{\kappa,d} < L^{\rm sc}_{\kappa,d}$ for all $d\geq2$ and $\kappa\geq3/2$. We deduce that in dimensions $2\leq d\leq 7$ the two continuous curves $L^{(1)}_{\kappa,d}$ and $L^{\rm sc}_{\kappa,d}$ must cross. The crossing point is unique by Lemma~\ref{al} and this concludes our proof of Theorem~\ref{thm:crossing}. 
\end{proof}


\section{Finite rank Lieb-Thirring inequalities: Proof of Theorem~\ref{thm:K}}
\label{sec:proof_K}

This section contains the proof of Theorem~\ref{thm:K} which, for convenience, we split into several intermediate steps. Our goal is to study the optimisation problem corresponding to inequality~\eqref{eq:LT_Schatten}, namely
\begin{equation}
\boxed{ K_{p,d}^{(N)} := \inf_{0 \le \gamma = \gamma^* \atop \rank(\gamma) \le N} \dfrac{ \norm{\gamma}_{\gS^q}^{\frac{p(2-d)+d}{d(p-1)}}\tr(-\Delta\gamma)}{\norm{\rho_\gamma}_{L^p(\R^d)}^{\frac{2p}{d(p-1)}}},}
\label{eq:def_K_N}
\end{equation}
where we recall that 
\begin{equation}
q := \begin{cases}
\frac{2p+d-dp}{2+d-dp}&\text{for $1<p<1+\frac2d$,}\\
+\ii&\text{for $p=1+\frac2d$.}
\end{cases}
\label{eq:relation_q}
\end{equation}
Throughout the paper, the constants $p$, $q$ and $\kappa$ are linked by the relations (we set $p' = p/(p-1)$ and $\kappa' = \kappa/ (\kappa - 1)$)
\[
    \boxed{\kappa +\frac{d}{2}= p' , \quad \text{and} \quad q = \kappa'.}
\]
Taking~\eqref{eq:def_K_N} to the power $\frac12(p-1)$, and letting $p \to 1$, so that $q \to 1$ as well, we recover the equality
\[
\int_{\R^d}\rho_\gamma(x)\,{\rd}x=\| \rho_\gamma \|_{L^1(\R^d)} = \| \gamma \|_{\gS^1}=\tr(\gamma),
\]
for all $0 \le \gamma = \gamma^*$. On the other hand, taking $p = 1 + 2/d$, so that $q = \infty$, we recover the better known dual Lieb-Thirring inequality
\begin{equation}
K_{1+2/d,d}^{(N)}\int_{\R^d}\rho_\gamma(x)^{1+\frac2d}\,\rd x\leq  \|\gamma\|^{\frac2d}\tr(-\Delta\gamma),\qquad \forall 0 \le \gamma=\gamma^*,\ \rank(\gamma)\leq N.
\label{eq:K_forkappa=1}
\end{equation}
We can think of~\eqref{eq:LT_Schatten} as a specific interpolation between these two cases. Note that a direct proof of~\eqref{eq:K_forkappa=1} with $N=+\ii$ can be found in~\cite{Rumin-11}, see also~\cite{LunSol-13,Sabin-16,Nam-18}. The original Lieb-Thirring proof proceeds by proving~\eqref{eq:LT_V_N} and then deducing~\eqref{eq:K_forkappa=1} by duality.

\subsection{Proof of $(i)$ on the existence of optimisers}
\label{sssec:step1}
Consider a minimising sequence $(\gamma_n)$ with $\rank(\gamma_n)\leq N$ for~\eqref{eq:def_K_N}, normalised such that
$$\tr(-\Delta\gamma_n)=1,\qquad \|\gamma_n\|_{\gS^q}=1$$
and
\begin{equation}
\lim_{n\to\ii}\int_{\R^d}\rho_{n}(x)^p\,\rd x=\frac1{\left(K^{(N)}_{p,d}\right)^{\frac{d(p-1)}{2}}}
\label{eq:min_seq_K}
\end{equation}
with $\rho_n:=\rho_{\gamma_n}$. We have $\|\gamma_n\|\leq \|\gamma_n\|_{\gS^q}=1$ and hence 
$$\int_{\R^d}\rho_n(x)\,\rd x=\tr(\gamma_n)\leq N.$$ 
This proves that $\rho_n$ is bounded in $L^1(\R^d)$. On the other hand, the Hoffmann-Ostenhof~\cite{Hof-77} inequality states that
\begin{equation}
\tr(-\Delta \gamma) \geq \int_{\R^d}|\nabla\sqrt{\rho_\gamma}(x)|^2\,\rd x
\label{eq:Hoffmann-Ostenhof}
\end{equation}
for all $\gamma=\gamma^*\geq0$. This shows that $\sqrt{\rho_n}$ is bounded in $H^1(\R^d)$, hence in $L^r(\R^d)$ for all $2\leq r<2^*$ where $2^*=2d/(d-2)$ in dimension $d\ge3$ and $2^*=+\ii$ in dimensions $d=1,2$, by the Sobolev inequality. In particular, we can choose $r = p$. From~\cite{Lieb-83} or from~\cite[Lem.~I.1]{Lions-84b}, we know that 
\begin{itemize}[leftmargin=*]
    \item \textbf{either} $\rho_n\to0$ strongly in $L^p(\R^d)$,
    \item \textbf{or} there is a $\rho\neq 0$ with $\sqrt\rho\in H^1(\R^d)$, a sequence $\tau_k\in\R^d$ and a subsequence so that $\sqrt{\rho_{n_k}(\cdot-\tau_k)}\wto \sqrt{\rho}\neq0$ weakly in $H^1(\R^d)$. 
\end{itemize}
Due to~\eqref{eq:min_seq_K} we know that the first possibility cannot happen and we may assume that $\sqrt\rho_n\wto \sqrt\rho\neq0$, after extraction of a subsequence and translation of the whole system by $\tau_n$. We may also extract a weak-$\ast$ limit for $\gamma_n$ in the trace class topology and infer $\gamma_n\wto\gamma$ where $\rho_\gamma=\rho\neq0$, hence $\gamma\neq0$. By passing to the limit, we have $\gamma=\gamma^*\geq0$ and $\rank(\gamma)\leq N$. 

Next we apply Lions' method~\cite{Lions-84} based on the Levy concentration function $Q_n(R)=\int_{|x|\leq R}\rho_n(x)\,\rd x$ and the strong local compactness in $L^2(\R^d)$ to deduce that there exists a sequence $R_n\to\ii$ so that
$$\lim_{n\to\ii}\int_{|x|\leq R_n}\rho_n(x)\,\rd x=\int_{\R^d}\rho(x)\,\rd x,\qquad \lim_{n\to\ii}\int_{R_n\leq |x|\leq 2R_n}\rho_n(x)\,\rd x=0.$$
Let $\chi\in C^\ii_c(\R^d,[0,1])$ be a smooth localisation function such that $\chi\equiv1$ on the unit ball $B_1$ and $\chi\equiv0$ outside of $B_2$. Let $\chi_n(x):=\chi(x/R_n)$ and $\eta_n=\sqrt{1-\chi_n^2}$. Then $\chi_n^2\rho_n\to\rho$ strongly in $L^1(\R^d)\cap L^p(\R^d)$ whereas $|\nabla\chi_n|^2\rho_n\to0$ and $|\nabla\eta_n|^2\rho_n\to0$
strongly in $L^1(\R^d)$. By the IMS formula (see, e.g., \cite[Thm.~3.2]{CycFroKirSim-87}) and Fatou's lemma for operators (see, e.g., \cite[Thm.~2.7]{Simon-79}), we obtain
\begin{align*}
\tr(-\Delta\gamma_n)&=\tr(-\Delta\chi_n\gamma_n\chi_n)+\tr(-\Delta\eta_n\gamma_n\eta_n)-\int_{\R^d}(|\nabla\chi_n|^2+|\nabla\eta_n|^2)\rho_n\\
&=\tr(-\Delta\chi_n\gamma_n\chi_n)+\tr(-\Delta\eta_n\gamma_n\eta_n)+o(1)\\
&\geq\tr(-\Delta\gamma)+\tr(-\Delta\eta_n\gamma_n\eta_n)+o(1).
\end{align*}
From the strong convergence of $\chi_n^2\rho_n$ we have
\begin{align*}
\int_{\R^d}\rho_n^p&=\int_{\R^d}\chi_n^2(\rho_n)^p+\int_{\R^d}(\eta_n^2\rho_n)^p+\int_{\R^d}(\eta_n^2-\eta_n^{2p})\rho_n^p\\
&=\int_{\R^d}\rho^p+\int_{\R^d}(\eta_n^2\rho_n)^p+o(1).
\end{align*}

First, we assume that $q<\ii$, that is, $p<1+2/d$. The Schatten norm satisfies 
\begin{align*}
\tr(\gamma_n)^q&=\tr\big(\chi_n(\gamma_n)^q\chi_n\big)+\tr\big(\eta_n(\gamma_n)^q\eta_n\big)\\
&\geq\tr(\chi_n\gamma_n\chi_n)^q+\tr(\eta_n\gamma_n\eta_n)^q\\
&\geq\tr(\gamma)^q+\tr(\eta_n\gamma_n\eta_n)^q+o(1).
\end{align*}
In the second line we have used the inequality $\tr(ABA)^m\leq \tr(A^mB^mA^m)$ for all $m\geq1$~\cite[App.~B]{LieThi-76} to infer 
$$\tr(\gamma_n)^q(\chi_n)^2\geq \tr(\gamma_n)^q(\chi_n)^{2q}=\tr(\chi_n)^q(\gamma_n)^q(\chi_n)^q\geq \tr(\chi_n\gamma_n\chi_n)^q.$$
In the third line we used Fatou's lemma in the Schatten space $\gS^q$.
Next, we argue using the method of the missing mass as in~\cite{Lieb-83c}, see also~\cite{Frank-13}, noticing that $K^{(N)}_{p,d}$ can be rewritten as
$$\left(K^{(N)}_{p,d}\right)^{\frac{d(p-1)}{2}}=\inf_{\substack{\gamma=\gamma^*\geq0\\ \rank(\gamma)\leq N}}\frac{\Big(\tr(\gamma^q)\Big)^{1-\theta}\Big(\tr(-\Delta\gamma)\Big)^{\theta}}{\int_{\R^d}\rho_\gamma(x)^{p}\,\rd x}$$
with
$$\theta:=\frac{d(p-1)}{2}\in(0,1).$$
Using H\"older's inequality in the form 
$$(a_1+a_2)^\theta(b_1+b_2)^{1-\theta}\geq a_1^\theta b_1^{1-\theta}+a_2^\theta b_2^{1-\theta}$$
we find 
\begin{align*}
1&=\Big(\tr(\gamma_n^q)\Big)^{1-\theta}\Big(\tr(-\Delta\gamma_n)\Big)^{\theta}\\
&\geq\Big(\tr(\gamma^q)\Big)^{1-\theta}\Big(\tr(-\Delta\gamma)\Big)^{\theta}+\Big(\tr(\eta_n\gamma_n \eta_n)^q\Big)^{1-\theta}\Big(\tr(-\Delta\eta_n\gamma_n\eta_n)\Big)^{\theta}+o(1)\\
&\geq\Big(\tr(\gamma^q)\Big)^{1-\theta}\Big(\tr(-\Delta\gamma)\Big)^{\theta}+\left(K^{(N)}_{p,d}\right)^{\frac{d(p-1)}{2}}\int_{\R^d}(\eta_n^2\rho_n)^p+o(1)\\
&=\Big(\tr(\gamma^q)\Big)^{1-\theta}\Big(\tr(-\Delta\gamma)\Big)^{\theta}+1-\left(K^{(N)}_{p,d}\right)^{\frac{d(p-1)}{2}}\int_{\R^d}\rho_{\gamma}^p+o(1).
\end{align*}
In the third line we used $\rank(\eta_n\gamma_n\eta_n)\leq N$. Passing to the limit we obtain
$$\left(K^{(N)}_{p,d}\right)^{\frac{d(p-1)}{2}}\int_{\R^d}\rho_{\gamma}^p\geq \Big(\tr(\gamma^q)\Big)^{1-\theta}\Big(\tr(-\Delta\gamma)\Big)^{\theta}$$
and therefore $\gamma\neq0$ is an optimiser. 

The case $p=1+2/d$ is similar. This time, we use $\|\gamma\|\leq \liminf_{n\to\ii}\|\gamma_n\|=1$ and $\|\eta_n\gamma_n\eta_n\|\leq\|\gamma_n\|=1$ to bound
\begin{align*}
1&=\tr(-\Delta\gamma_n)\\
&\geq\tr(-\Delta\gamma)+\tr(-\Delta\eta_n\gamma_n\eta_n)+o(1)\\
&\geq\|\gamma\|^{\frac2d}\tr(-\Delta\gamma)+\|\eta_n\gamma_n\eta_n\|^{\frac2d}\tr(-\Delta\eta_n\gamma_n\eta_n)+o(1)\\
&\geq\|\gamma\|^{\frac2d}\tr(-\Delta\gamma)+K_{1+2/d,d}^{(N)}\int_{\R^d}(\eta_n^2\rho_n)^{1+\frac2d}+o(1)\\
&=\|\gamma\|^{\frac2d}\tr(-\Delta\gamma)+1-K_{1+2/d,d}^{(N)}\int_{\R^d}\rho_\gamma^{1+\frac2d}+o(1)
\end{align*}
and arrive at the same conclusion that $\gamma$ is an optimiser.

\subsection{Proof of $(ii)$ on the equation}

Let $\gamma$ be an optimiser such that 
$$\tr(-\Delta\gamma)=\int_{\R^d}\rho(x)^p\,\rd x=1.$$
This normalisation is always possible by scaling and by multiplying $\gamma$ by a positive constant. Then we have
$$\tr(\gamma^q)=\left(K^{(N)}_{p,d}\right)^{\frac{d(p-1)}{2+d-dp}}.$$

We start with the case $q<\ii$, that is, $p<1+2/d$. Assume that we have a smooth curve of operators $\gamma(t)=\gamma+t\delta+o(t)$ for some $\delta=\delta^*$, with $\gamma(t)=\gamma(t)^*\geq0$ and $\rank(\gamma(t))\leq N$. By expanding we find
\begin{align}
\left(K^{(N)}_{p,d}\right)^{\frac{d(p-1)}{2}}&\leq \frac{\Big(\tr(\gamma(t)^q)\Big)^{1-\theta}\Big(\tr(-\Delta\gamma(t))\Big)^{\theta}}{ \int_{\R^d}\rho_{\gamma(t)}^p}\nn\\
&=\left(K^{(N)}_{p,d}\right)^{\frac{d(p-1)}{2}}\frac{\Big(1+qt\frac{\tr(\delta\gamma^{q-1})}{\tr(\gamma^q)}+o(t)\Big)^{1-\theta}\Big(1+t\tr(-\Delta\delta)+o(t)\Big)^{\theta}}{1+pt\int_{\R^d}\rho_\delta\rho_\gamma^{p-1}+o(t) }\nn\\
&=\left(K^{(N)}_{p,d}\right)^{\frac{d(p-1)}{2}}\!\!\left(1+t\,\theta\,\tr\left[\delta\left(-\Delta-\frac{p}{\theta}\rho_\gamma^{p-1}+\frac{q(1-\theta)}{\theta\tr(\gamma^q)}\gamma^{q-1}\right)\right]+o(t)\right).\label{eq:derivative}
\end{align}
Now take $\gamma(t):=e^{itH}\gamma e^{-itH}=\gamma+it[H,\gamma]+o(t)$ for some (smooth and finite-rank) self-adjoint operator $H$ and all $t\in\R$. Since $\rank(\gamma(t))=\rank(\gamma)$, we deduce from~\eqref{eq:derivative} after varying over all $H$ that 
$$\left[-\Delta-\frac{p}{\theta}\rho_\gamma^{p-1}\,,\,\gamma\right]=0.$$
Hence $\gamma$ commutes with the mean-field operator $H_\gamma:=-\Delta-p\rho_\gamma^{p-1}/\theta$. We can therefore write $\gamma=\sum_{j=1}^R n_j|u_{k_j}\rangle\langle u_{k_j}|$ for some eigenvectors $u_{k_j}$ of $H_\gamma$ (with eigenvalue $\mu_{k_j}$) and some $n_j>0$. In particular, $H_\gamma$ admits at least $R$ eigenvalues. 

Using now $\gamma(t)=\gamma+t\delta$ for a $\delta$ supported on the range of $\gamma$ and for $t$ small enough in~\eqref{eq:derivative}, we find that 
$$-\Delta-\frac{p}{\theta}\rho_\gamma^{p-1}+\frac{(1-\theta)q}{\theta\tr(\gamma^q)}\gamma^{q-1}\equiv0\qquad\text{on the range of $\gamma$.}$$
Evaluating this identity on $u_{k_j}$ we infer that
$$
\mu_{k_j} + \frac{(1-\theta)q}{\theta\tr(\gamma^q)} n_j^{q-1} =0.
$$
This shows that $\mu_{k_j}<0$ and 
$$n_j = \left(\frac{\theta\tr(\gamma^q)}{(1-\theta)q}\right)^{\frac1{q-1}}\ | \mu_{k_j} |^{\frac{1}{q - 1}}.$$
Since $\gamma$ is assumed to be of rank $R$, we in particular deduce that $H_\gamma$ has at least $R$ negative eigenvalues.
 
Next, we show that the $\mu_{k_j}$ are necessarily the $R$ first eigenvalues. Assume that one eigenvector of $H_\gamma$ with eigenvalue $<\mu_{R}$ does not belong to the range of $\gamma$, so there is $1 \le j \le R$ with $u_{k_j} \neq u_j$ with $k_j > j$ and $u_j$ not in the range of $\gamma$. Consider the new operator
$$\gamma':=\gamma - n_j |u_{k_j}\rangle\langle u_{k_j}|+ n_j |u_j\rangle\langle u_j|:=\gamma+\delta,$$
which has the same rank and the same $\gS^q$ norm as $\gamma$. We have by convexity 
$$\int_{\R^d}\rho_{\gamma'}^p\geq 1+pn_{j}\int_{\R^d}\rho_\gamma^{p-1}\left(|u_j|^2-|u_{k_j}|^2\right)$$
and 
\begin{align*}
\tr(-\Delta\gamma')&=1+n_{j}\pscal{u_j,-\Delta u_j}-n_{k_j}\pscal{u_{k_j},-\Delta u_{k_j}} \\
&=1+\frac{pn_{j}}{\theta}\int_{\R^d}\rho^{p-1}_\gamma \big(|u_j|^2-|u_{k_j}|^2\big)+\left(\mu_j-\mu_{k_j}\right) n_j \\
&<1+\frac{pn_{j}}{\theta}\int_{\R^d}\rho^{p-1}_\gamma \big(|u_j|^2-|u_{k_j}|^2\big)
\end{align*}
since $\mu_j<\mu_{k_j}$. This gives
\begin{align*}
\frac{\Big(\tr(\gamma')^q\Big)^{1-\theta}\Big(\tr(-\Delta\gamma')\Big)^\theta}{\int_{\R^d}\rho_{\gamma'}^p}
&<\left(K^{(N)}_{p,d}\right)^{\frac{d(p-1)}{2}}\frac{\left(1+\frac{pn_{j}}{\theta}\int_{\R^d}\rho^{p-1}_\gamma \big(|u_j|^2-|u_{k_j}|^2\big)\right)^\theta}{1+pn_{j}\int_{\R^d}\rho_\gamma^{p-1}\left(|u_j|^2-|u_{k_j}|^2\right)} \\
&\leq \left(K^{(N)}_{p,d}\right)^{\frac{d(p-1)}{2}},
\end{align*}
a contradiction. Hence $\mu_{k_j}=\mu_j$.

Finally, when $R<N$ and $\mu_{R+1}<0$, we can consider the operator 
$$\gamma(t)=\gamma+t|u_{R+1}\rangle\langle u_{R+1}|$$
with $t\geq 0$, which has rank $R+1\leq N$. From~\eqref{eq:derivative} we obtain
\begin{align*}
\left(K^{(N)}_{p,d}\right)^{\frac{d(p-1)}{2}}&\leq \left(K^{(N)}_{p,d}\right)^{\frac{d(p-1)}{2}}\bigg(1+o(t)\\
&\qquad +t\theta\pscal{u_{R+1},\left(-\Delta-\frac{p}{\theta}\rho_\gamma^{p-1}+\frac{(1-\theta)q}{\theta\tr(\gamma^q)}\gamma^{q-1}\right)u_{R+1}}\bigg)\\
&\leq \left(K^{(N)}_{p,d}\right)^{\frac{d(p-1)}{2}}\left(1+t\mu_{R+1}\theta+o(t)\right),
\end{align*}
another contradiction. Hence $H_\gamma$ cannot have more than $R$ negative eigenvalues when $R<N$.

As a conclusion, we have shown that 
$$\gamma=\left(\frac{\theta\tr(\gamma^q)}{q(1-\theta)}\right)^{\frac{1}{q-1}}\sum_{j=1}^R|\mu_j|^{\frac{1}{q-1}}|u_j\rangle\langle u_j|,$$
with
$$\left(-\Delta-\frac{p}{\theta}\rho_\gamma(x)^{p-1}\right)u_j=\mu_j\,u_j,\qquad j=1,...,R.$$
Taking the trace of $\gamma^q$ we find that 
$$\frac{\theta\tr(\gamma^q)}{q(1-\theta)}=\left(\frac{q(1-\theta)}{\theta} \dfrac{1}{\sum_{j=1}^R|\mu_j|^{\frac{q}{q-1}} }\right)^{q-1}$$
and thus
$$\gamma=\frac{q(1-\theta)}{\theta\sum_{j=1}^R|\mu_j|^{\frac{q}{q-1}}} \sum_{j=1}^R|\mu_j|^{\frac{1}{q-1}}|u_j\rangle\langle u_j|.$$
Replacing $\gamma$ by $(p/\theta)^{\frac1{p-1}}\gamma$ we find the equation mentioned in the statement.

\medskip

The arguments for $q=+\ii$ ($p=1+2/d$) are similar. We start with a minimiser normalised so that 
$$\int_{\R^d}\rho_\gamma^{1+\frac2d}=\tr(-\Delta\gamma)=1,\qquad \|\gamma\|^{\frac2d}=K_{1+2/d,d}^{(N)}.$$ 
The first perturbation $\gamma(t):=e^{itH}\gamma e^{-itH}=\gamma+it[H,\gamma]+o(t)$ leaves the operator norm invariant and provides the equation $[-\Delta-p\rho_\gamma^{2/d}\,,\,\gamma]=0$, hence again $\gamma=\sum_{j=1}^Rn_j|u_{k_j}\rangle\langle u_{k_j}|$ with $H_\gamma u_{k_j}=\mu_{k_j}u_{k_j}$ and $H_\gamma=-\Delta-p\rho_\gamma^{2/d}$. 
In order to prove that $\mu_{k_j}<0$, we consider the operator
$$\tilde\gamma:=\gamma-n_j|u_{k_j}\rangle\langle u_{k_j}|$$
which has one less eigenvalue and satisfies $\|\tilde\gamma\|^{2/d}\leq \|\gamma\|^{2/d}=K^{(N)}_{1+2/d,d}$. We find
\begin{align*}
K^{(N)}_{1+2/d,d}\leq K^{(N-1)}_{1+2/d,d}&\leq \frac{\|\tilde \gamma\|^{\frac2d}\tr(-\Delta\tilde \gamma)}{ \int_{\R^d}\rho_{\tilde\gamma}^{1+\frac2d}}\nn\\
&\leq K^{(N)}_{1+2/d,d}\frac{\tr(-\Delta\tilde \gamma)}{ \int_{\R^d}\rho_{\tilde\gamma}^{1+\frac2d}}\nn\\
&=K^{(N)}_{1+2/d,d}\frac{1-n_j\int_{\R^d}|\nabla u_{k_j}|^2}{\int_{\R^d}\big(\rho_\gamma-n_j|u_{k_j}|^2)^{1+\frac2d}}\nn\\
&=K^{(N)}_{1+2/d,d}\frac{1-n_j\mu_{k_j}-n_j\frac{d+2}{d}\int_{\R^d}\rho_\gamma^{\frac2d}|u_{k_j}|^2}{\int_{\R^d}\big(\rho_\gamma-n_j|u_{k_j}|^2)^{1+\frac2d}}.\label{eq:decrease_occ}
\end{align*}
Simplifying by $K^{(N)}_{1+2/d,d}>0$, this gives the estimate
\begin{equation}
\mu_{k_j}\leq -\frac1{n_j}\left(\int_{\R^d}\big(\rho_\gamma-n_j|u_{k_j}|^2)^{1+\frac2d}-\int_{\R^d}\rho_{\gamma}^{1+\frac2d}+n_j\frac{d+2}{d}\int_{\R^d}\rho_\gamma^{\frac2d}|u_{k_j}|^2\right)<0
\label{eq:estim_mu_j}
\end{equation}
where the last negative sign is by strict convexity of $t\mapsto t^{1+2/d}$. 
Hence $\gamma$ has its range into the negative spectral subspace of $H_\gamma$, an operator which thus possesses at least $R$ negative eigenvalues. 
Next we show that $n_j=\|\gamma\|$ for all $j=1,...,R$. Assume on the contrary that $0<n_j<\|\gamma\|$ (this can only happen when $R\geq2$). Taking $\gamma(t)=\gamma+t|u_{k_j}\rangle\langle u_{k_j}|$ which has the same operator norm for $t$ small enough, we obtain 
\begin{align}
K^{(N)}_{1+2/d,d}\leq \frac{\|\gamma(t)\|^{\frac2d}\tr(-\Delta\gamma(t))}{ \int_{\R^d}\rho_{\gamma(t)}^{1+\frac2d}}&=K^{(N)}_{1+2/d,d}\frac{1+t\int_{\R^d}|\nabla u_{k_j}|^2}{\int_{\R^d}\big(\rho_\gamma+t|u_{k_j}|^2)^{1+\frac2d}}\nn\\
&=K^{(N)}_{1+2/d,d}\frac{1+t\mu_{k_j}+pt\int_{\R^d}\rho_\gamma^{p-1}|u_{k_j}|^2}{\int_{\R^d}\big(\rho_\gamma+t|u_{k_j}|^2)^{1+\frac2d}}\nn\\
&=K^{(N)}_{1+2/d,d}\left(1+t\mu_{k_j}+o(t)\right)\label{eq:increase_occ}
\end{align}
which is a contradiction since $\mu_{k_j}<0$, as we have seen. We conclude that $n_j=\|\gamma\|$ for all $j=1,...,R$. The argument for showing that $\mu_{k_1},...,\mu_{k_R}$ are the $R$ first eigenvalues is exactly the same as before.

\subsection{Proof of $(iii)$ on the rank of optimisers}

In this subsection, we prove the following result.

\begin{proposition}[Binding]\label{prop:binding}
    Let $1<p\leq 1+2/d$ with $p<2$ and assume that $K^{(N)}_{p,d}$ admits an optimiser $\gamma$ of rank $N$. Then 
    $K^{(2N)}_{p,d}<K^{(N)}_{p,d}$.
\end{proposition}

The proof of $(iii)$ in Theorem~\ref{thm:K} follows immediately from Proposition~\ref{prop:binding}, arguing as follows. Since $K^{(1)}_{p,d}$ has an optimiser, the proposition shows that $K^{(2)}_{p,d}<K^{(1)}_{p,d}$, hence we can take $N_2=2$. By Step $(i)$ there is an optimiser for $K^{(2)}_{p,d}$ and by Step $(ii)$ the strict inequality $K^{(2)}_{p,d}<K^{(1)}_{p,d}$ implies that the optimisers for $K^{(2)}_{p,d}$ all have rank two. Hence Proposition~\ref{prop:binding} implies that $K^{(4)}_{p,d}<K^{(2)}_{p,d}$. If $K^{(3)}_{p,d}<K^{(2)}_{p,d}$ we take $N_3=3$ and otherwise we take $N_3=4$. We then go on by induction to obtain the assertion of $(iii)$. Hence we now concentrate on proving Proposition~\ref{prop:binding}. 

\begin{proof}[Proof of Proposition~\ref{prop:binding}]
    We follow ideas from~\cite[Section~2.4]{GonLewNaz-20_ppt}. Let $\gamma := \sum_{j=1}^{N}n_j | u_j \ket \bra u_j |$  be a minimiser of rank $N$ for $K^{(N)}_{p,d}$, normalised in the manner $\tr(-\Delta\gamma)=\int_{\R^d}\rho^p=1$. The functions $u_j$ satisfy 
    $$\left(-\Delta-\frac{p}{\theta}\left(\sum_{j=1}^Nn_j|u_j|^2\right)^{p-1}\right)u_j=\mu_j\,u_j$$
    with $n_j=c|\mu_j|^{1/(q-1)}$. Note that the first eigenfunction $u_1$ is positive, hence the nonlinear potential never vanishes. By usual regularity arguments, this shows that the $u_j$ are $C^\ii$ and decay exponentially at infinity. For $R>0$, we set $u_{j,R}(x) := u_j(x - R e_1)$ where $e_1=(1,0,...,0)$, and we introduce the Gram matrix
    \[
    S_R = \begin{pmatrix}
    \bbI_N & E^R \\
    (E^{R})^* & \bbI_N
    \end{pmatrix}, \quad \text{with} \quad E_{ij}^R := \bra u_i, u_{j,R} \ket = \int_{\R^d} u_i(x) u_j(x - Re_1) \rd x.
    \]
    Since the functions $u_i$ are exponentially decaying, $E_R$ goes to $0$, and the overlap matrix $S_R$ is invertible for $R$ large enough. We then let 
    $$\begin{pmatrix}
    \psi_{1,R}\\ \vdots\\ \psi_{2N,R}
    \end{pmatrix}=(S_R)^{-\frac12}\begin{pmatrix}u_1\\ \vdots\\ u_N\\ u_{1,R}\\ \vdots\\ u_{N,R}\end{pmatrix}$$
    and
    $$\gamma_R = \sum_{j=1}^{N} n_j\Big(| \psi_{j,R} \ket \bra \psi_{j,R} |+| \psi_{N+j,R} \ket \bra \psi_{N+j,R} |\Big).$$
    We have
    $$\tr(\gamma_R)^q=2\tr(\gamma^q),\qquad \|\gamma_R\|=\|\gamma\|.$$
    Expanding as in~\cite{GonLewNaz-20_ppt} using 
    \[
    (S_R)^{-1/2} = \begin{pmatrix}
    \bbI_N & 0 \\
    0 & \bbI_N
    \end{pmatrix} - \frac12\begin{pmatrix}
    0 & E^R \\
    (E^R)^* & 0
    \end{pmatrix} + \frac38\begin{pmatrix}
    E^R (E^R)^* & 0 \\
    0 & (E^R)^* E^R
    \end{pmatrix} + O(e_R^3).
    \]
    for
    $$e_R:=\max_{i,j}\int_{\R^d} |u_i(x)|\,|u_j(x - Re_1)| \rd x,$$     
    we obtain after a long calculation
    \begin{align*}
    \left(K^{(2N)}_{p,d}\right)^{\frac{d(p-1)}{2}} &\leq \left(K^{(N)}_{p,d}\right)^{\frac{d(p-1)}{2}} \frac{2^{1-\theta}\big(\tr(-\Delta\gamma_R)\big)^\theta}{\int_{\R^d}\rho_{\gamma_R}^p}\\
    &=\left(K^{(N)}_{p,d}\right)^{\frac{d(p-1)}{2}} \left(1-\frac12\int_{\R^d} \left( (\rho+\rho_R)^p-\rho^p-\rho_R^p \right) +O(e_R^2)\right)
    \end{align*}
    with $\rho(x)=\rho_\gamma(x)$ and $\rho_R(x)=\rho(x-Re_1)$. 
    From the arguments in~\cite[Section~2.4]{GonLewNaz-20_ppt} we know that 
    \begin{equation}
    \label{eq:boundgln}
    \int_{\R^d} \left( (\rho+\rho_R)^p-\rho^p-\rho_R^p \right) \geq c R^{p(1-d)}e^{-p\sqrt{|\mu_N|}R}
    \end{equation}
    and by~\cite[Lemma~21]{GonLewNaz-20_ppt} we have 
    $$e_R\leq C(1+R^d)e^{-\sqrt{|\mu_N|}R}.$$
    Since $p<2$ by assumption we conclude, as we wanted, that $K^{(2N)}_{p,d}<K^{(N)}_{p,d}$.
\end{proof}

\section{Binding for $\kappa<1$ and $N=2$: Proof of Theorem~\ref{thm:LT_bis}}
\label{sec:proof_LT_V_kappa<1}

In this section we provide the proof of Theorem~\ref{thm:LT_bis}. Define $p$ by $p'=\kappa+d/2$ let $Q$ be the radial Gagliardo--Nirenberg minimiser, solution to~\eqref{eq:NLS}, and set $m:=\int_{\R^d} Q^2\,dx$. 

\subsection{Some properties of $Q$}
First we relate our constants for $N = 1$ to $Q$. We have the Pohozaev identity
\begin{equation}
\label{eq:pohozaev}
\begin{cases}
\dps \int_{\R^d} |\nabla Q|^2\,dx - \int_{\R^d} Q^{2p}\,dx = -m,\\[0.4cm]
\dps \left( \frac d2-1 \right)\int_{\R^d} |\nabla Q|^2\,dx - \frac d{2p} \int_{\R^d} Q^{2p}\,dx = - \frac d2 m \,. 
\end{cases}
\end{equation}
These follow by multiplying the equation~\eqref{eq:NLS} by $ Q$ and by $x\cdot\nabla Q$, respectively. This gives the identity
\begin{equation} \label{eq:intQ2p}
    \frac{m}{\int_{\R^d} Q^{2p}} = 1 - \frac{d}{2} \frac{p-1}{p} = \frac{p-1}{p} \kappa.
\end{equation}
On the other hand, setting $V_Q := -Q^{2(p-1)}$, we see that $Q$ is an eigenvector of $-\Delta + V_Q$ (with corresponding eigenvalue $-1$), and, by optimality of $V_Q$ for $L^{(1)}_{\kappa, d}$, we have
\begin{equation} \label{eq:exactL1_withQ}
    L^{(1)}_{\kappa, d} = \frac{1}{\int_{\R^d} | V_Q |^{\kappa + \frac{d}{2}}} = \frac{1}{\int_{\R^d} Q^{2p}}.
\end{equation}
Finally, it is well known that there is $C > 0$ so that
\begin{equation} \label{eq:boundQ}
\frac{1}{C} \frac{\re^{ - | x |}}{1 + | x |^{\frac{d-1}{2}}} \le Q(x) \le C \frac{\re^{ - | x |}}{1 + | x |^{\frac{d-1}{2}}}.
\end{equation}

\subsection{Test potential for $L^{(2)}_{\kappa, d}$}

We now construct a test potential to find a lower bound for $L^{(2)}_{\kappa, d}$.  For $R>0$, We let
$$
Q_\pm(x) = Q\big(x\pm\tfrac R2 e_1\big) 
$$
with $e_1=(1,0,...,0)$. Inspired by the dual problem studied in the previous section, we consider the potential
$$
\boxed{V = - \left(Q_+^2 + Q_-^2\right)^{p-1} \,.}
$$
It is important here that we add the two densities and not the corresponding potentials. We do not see how to make our proof work if we would take $V = -Q_+^{2(p-1)} - Q_-^{2(p-1)}$ instead. 

We introduce the quantity
\begin{equation}
    \label{eq:def:A}
    A = A(R) := \frac12 \int_{\R^d} \left( (Q_+^2+Q_-^2)^p - Q_+^{2p} - Q_-^{2p} \right)dx >0\,.	
\end{equation}
Due to the inequality~\eqref{eq:boundQ}, $A$ goes (exponentially fast) to $0$ as $R$ goes to infinity. Our main result is the following.
\begin{lemma}
    We have, as $R \to \infty$,
    \[
    L^{(2)}_{\kappa, d} \ge \frac{|\lambda_1(-\Delta+V)|^\kappa + |\lambda_2(-\Delta+V)|^\kappa}{\int_{\R^d} |V |^{\kappa+\frac{d}{2}}\,dx} 
    = L^{(1)}_{\kappa, d} \left(1 + \frac{\kappa}{pm} A + o(A) \right).
    \]
\end{lemma}
The proof of Theorem~\ref{thm:LT_bis} follows as the leading correction is positive.

\begin{proof}
	First, we bound $A$ from below similarly to \eqref{eq:boundgln}. Indeed, noting that the integrand of $A$ is nonnegative and bounding it from below using~\eqref{eq:boundQ} in a neighborhood of the origin, we find
\begin{equation} \label{eq:boundA}
A \ge \frac12 \int_{\cB(0, 1)}  \left( (Q_+^2+Q_-^2)^p - Q_+^{2p} - Q_-^{2p} \right) \ge c \dfrac{\re^{ - p R}}{R^{p(d-1)}}.
\end{equation}

Next, we turn to the denominator appearing in the lemma. We have
$$
\int_{\R^d} | V |^{\kappa+\frac{d}{2}}\,dx = \int_{\R^d} \left(Q_+^2 + Q_-^2\right)^{p}  = 2 \int_{\R^d} Q^{2p}\,dx  + 2A.
$$
Together with~\eqref{eq:exactL1_withQ}, this gives
\begin{align*}
\dfrac{1}{\int_{\R^d} | V |^{\kappa+\frac{d}{2}}\,dx} & = \frac{1}{2} \frac{1}{\int_{\R^d} Q^{2p}} \left( 1 - \frac{A}{\int_{\R^d} Q^{2p}} + O(A^2) \right) \\
& = \frac{L^{(1)}_{\kappa, d}}{2}  \left( 1 - \frac{A}{\int_{\R^d} Q^{2p}} + O(A^2) \right).
\end{align*}

Finally, we evaluate the numerator. We set $E := E(R) = \int_{\R^d} Q_+ Q_-\,dx$. We have $E \to 0$ as $R\to\infty$, so for $R$ large enough, we have $|E|<m$, and the two functions $\psi^{(\pm)}$ defined by
$$
\begin{pmatrix}
\psi^{(+)} \\ \psi^{(-)}
\end{pmatrix}
= \begin{pmatrix}
m & E \\ E & m
\end{pmatrix}^{-1/2} 
\begin{pmatrix}
Q_+ \\ Q_- 
\end{pmatrix}
$$
are orthonormal in $L^2(\R^d)$.	Let
$$
\mathcal H := \begin{pmatrix}
\langle \psi^{(+)},(-\Delta+V)\psi^{(+)}\rangle & \langle \psi^{(+)},(-\Delta+V)\psi^{(-)}\rangle \\
\langle \psi^{(-)},(-\Delta+V)\psi^{(+)}\rangle & \langle \psi^{(-)},(-\Delta+V)\psi^{(-)}\rangle
\end{pmatrix} \,.
$$
By the variational principle, the two lowest eigenvalues of $-\Delta+V$ are not larger than the corresponding eigenvalues of $\mathcal H$, and therefore
$$
|\lambda_1(-\Delta+V)|^\kappa + |\lambda_2 (-\Delta+V)|^\kappa \geq \Tr \;\mathcal H_-^\kappa \,.
$$
We have
$$
\mathcal H = h \bbI_2 + \begin{pmatrix}
0 & \delta \\ \delta & 0
\end{pmatrix},
$$
where
$$
h := \langle \psi^{(+)},(-\Delta+V)\psi^{(+)}\rangle = \langle \psi^{(-)},(-\Delta+V)\psi^{(-)}\rangle
$$
and
$$
\delta := \langle \psi^{(+)},(-\Delta+V)\psi^{(-)}\rangle = \langle \psi^{(-)},(-\Delta+V)\psi^{(+)}\rangle \,.
$$
We have $h\to -1$ and $\delta\to 0$ as $R\to\infty$, and therefore
$$
\Tr\; \mathcal H_-^\kappa = 2|h|^{\kappa} - \kappa |h|^{\kappa-1} \Tr  \begin{pmatrix}
0 & \delta \\ \delta & 0
\end{pmatrix} 
+ O(\delta^2) 
= 2|h|^{\kappa} + O(\delta^2) \,.
$$

It remains to expand $h$ and to bound $\delta$. We begin with $h$. We find
\begin{align*}
|\nabla\psi^{(+)}|^2 + |\nabla\psi^{(-)}|^2 & = \frac{m}{m^2-E^2} \left( |\nabla Q_+|^2 + |\nabla Q_-|^2 \right)  - \frac{2E}{M^2-E^2} \nabla Q_+ \cdot \nabla Q_-.
\end{align*}
Integrating and using~\eqref{eq:NLS} gives
\begin{align*}
\int_{\R^d} \left( |\nabla\psi^{(+)}|^2 + |\nabla\psi^{(-)}|^2 \right)dx & = -2 + \frac{2m}{m^2-E^2} \int_{\R^d} Q^{2p}\,dx \\
& \ \quad - \frac{E}{m^2-E^2} \int_{\R^d} \left( Q_+^{2p-2}+Q_-^{2p-2}\right)Q_+ Q_- \,dx \,.
\end{align*}
Similarly,
$$
(\psi^{(+)})^2 + (\psi^{(-)})^2 = \frac{m}{m^2-E^2} \left( Q_+^2 + Q_-^2 \right)  - \frac{2E}{M^2-E^2} Q_+ Q_-
$$
and therefore
\begin{align*}
h & = \frac12 \left( \langle \psi^{(+)},(-\Delta+V)\psi^{(+)}\rangle + \langle \psi^{(-)},(-\Delta+V)\psi^{(-)}\rangle \right) \\
& = -1 - \frac{m}{m^2-E^2} A + \frac{E}{m^2-E^2} B \,,
\end{align*}
where $A$ was defined in~\eqref{eq:def:A}, and where
$$
B=B(R) := \int_{\R^d} Q_+Q_- \left( (Q_+^2+Q_-^2)^{p-1} - \frac12 \left( Q_+^{2p-2} + Q_-^{2p-2}\right)\right)dx \,.
$$

From~\eqref{eq:boundQ} and~\cite[Lem.~21]{GonLewNaz-20_ppt} we see that $E(R) \le C' R^{d}\re^{-R}$ and $B(R) \le C' R^{d}\re^{ - R}$. In particular, by \eqref{eq:boundA} and the assumption $p < 2$, we have $E^2 = o(A)$ and $E B = o(A)$. This gives
\begin{align*}
    |h|^\kappa = (-h)^\kappa &= (1+ m^{-1}A + o(A))^\kappa = 1 + \kappa m^{-1} A + o(A) \,.
\end{align*}
We see in a similar fashion that $\delta \le C'R^d \re^{-R}$ hence $O(\delta^2) = o(A)$ as well. Gathering all the estimates gives
\[
    L^{(2)}_{\kappa, d} \ge L^{(1)}_{\kappa, d}
    \left(1 + \left( \kappa - \frac{m}{\int_{\R^d} Q^{2p}} \right) \frac{A}{m} + o(A) \right)
    =  L^{(1)}_{\kappa, d} \left(1 +\frac{\kappa}{pm} A + o(A) \right),
\]
where the last equality comes from~\eqref{eq:intQ2p}.    
\end{proof}

\section{Non existence of minimisers for the Fermionic NLS: Proof of Theorems~\ref{th:nonExistence_J(N)} and~\ref{th:N=2}}
\label{sec:NLS}

In this section, we prove our results concerning the minimisation problem $J(N)$ which, we recall, is defined by
\begin{equation}
J(N):=\inf \Big\{ \tr(-\Delta\gamma)-\frac1p\int_{\R^d}\rho_\gamma(x)^p\,\rd x : \ 0\leq \gamma=\gamma^*\leq1,\ \Tr(\gamma) = N \Big\}.
\label{eq:def_J}
\end{equation}
We assume in the whole section
\[
    1 < p < 1 + \frac2d.
\]
After an appropriate scaling, and using the fact that $\Tr(\gamma) = \| \gamma \|_{\gS^1}$, the optimal inequality $\cE(\gamma)\geq J(N)$ becomes
\begin{equation*}
\widetilde K_{p,d}^{(N)} \| \rho_\gamma \|_p^{\frac{2p}{d(p-1)}} \leq \| \gamma \|_{\gS^1}^{\frac{d+2-dp}{d(p-1)}}\; \tr(-\Delta\gamma),
\end{equation*}
valid for all $0 \le \gamma = \gamma^* \le 1$ with $\Tr(\gamma) = N$, and with best constant
\begin{equation} \label{eq:explicit_tildeK}
\boxed{\widetilde{K}_{p,d}^{(N)} := \left(\frac{|J(N)|}{N}\right)^{-\frac{d+2-pd}{d(p-1)}}\frac1{p-1}\left(\frac{d}{2p}\right)^{\frac{2}{d(p-1)}} \left(1+\frac2d-p\right)^{-\frac{d+2-dp}{d(p-1)}}.}
\end{equation}
One can remove the constraint $\| \gamma \| \le 1$ at the expense of a factor $\| \gamma \|^{d/2}$, and we obtain the optimal inequality
\begin{equation}
\boxed{ \widetilde K_{p,d}^{(N)} \| \rho_\gamma \|_p^{\frac{2p}{d(p-1)}} \leq \| \gamma \|_{\gS^1}^{\frac{d+2-dp}{d(p-1)}}\;\|\gamma\|^{\frac{2}{d}}\;\tr(-\Delta\gamma),}
\label{eq:LT_sub_critical}
\end{equation}
valid for all $0 \le \gamma = \gamma^*$ with $\Tr(\gamma) = N$.

\subsection{Link between NLS and Lieb-Thirring, proof of Theorem~\ref{th:nonExistence_J(N)}}
\label{ssec:proof_nonExistence_J(N)}

The link between the constant $\widetilde K_{p,d}^{(N)}$ and the dual Lieb-Thirring constant $K_{p, d}^{(N)}$ defined in~\eqref{eq:LT_Schatten} is given in the following proposition.

\begin{proposition}[Relation between $\widetilde K_{p,d}^{(N)}$ and $K_{p,d}^{(N)}$]\label{prop:relation_K}
    Let $d\geq1$ and $1<p<1+\frac2d$. For all $N\in\N$ we have 
    \begin{equation}
    K_{p,d}^{(N)}\leq \widetilde K_{p,d}^{(N)} \le \widetilde K_{p,d}^{(1)}=K_{p,d}^{(1)}.
    \label{eq:relation_K}
    \end{equation}
\end{proposition}

\begin{proof}
    It is shown in~\cite[Lemma~11]{GonLewNaz-20_ppt} that the minimisation problem $J(N)$ can be restricted to operators $\gamma$ which are orthogonal projectors of rank $N$. For such operators, we have $\| \gamma \| = 1$ and
    \[
        \| \gamma \|_{\gS^q}^q = \tr(\gamma^q) = N = \| \gamma \|_{\gS^1} = \rank(\gamma).
    \]
This gives
    \begin{equation*}
    K^{(N)}_{p,d}\leq \frac{\norm{\gamma}_{\gS^q}^{\frac{p(2-d)+d}{d(p-1)}}\tr(-\Delta\gamma)}{\norm{\rho_\gamma}_{L^p(\R^d)}^{\frac{2p}{d(p-1)}}} = \frac{ \| \gamma \|_{\gS^1}^{\frac{d+2-dp}{d(p-1)}}\|\gamma\|^{\frac2d}\tr(-\Delta\gamma)}{\norm{\rho_\gamma}_{L^p(\R^d)}^{\frac{2p}{d(p-1)}}}.
    \end{equation*}
    Optimising over projectors $\gamma$ gives $K_{p,d}^{(N)}\leq \widetilde K_{p,d}^{(N)}$. In the case $N = 1$, every operator of rank $1$ is proportional to a rank $1$ projector, so the two problems coincide, and $\widetilde{K}^{(1)}_{p,d}=K^{(1)}_{p,d}$. Finally, in~\cite{GonLewNaz-20_ppt}, it is also proved that $J(N) \le N J(1)$. This implies $\widetilde{K}_{p,d}^{(N)} \le \widetilde K_{p,d}^{(1)}$.
\end{proof}

There is a similarity between the proof of the above proposition and the arguments in \cite{LieLla-78,FraLieSeiTho-11b}. In those works also the sharp Lieb-Thirring inequality for $\kappa=3/2$ is used to obtain an inequality about orthonormal functions.

The relation~\eqref{eq:relation_K} allows us to prove Theorem~\ref{th:nonExistence_J(N)}, which states that $J(N) = N J(1)$ for all $N \in \N$, and that $J(N)$ admits no minimiser for $N \ge 2$.

\begin{proof}[Proof of Theorem~\ref{th:nonExistence_J(N)}]
It was proved in~\cite{LieThi-76} that for $\kappa = 3/2$, we have $L_{3/2,1}=L^{(N)}_{3/2,1}=L^{(1)}_{3/2,1}$ for all $N\in\N$. This implies $K_{2,1}^{(N)}=K_{2,1}^{(1)}$ for all $N \in \N$. Hence, by~\eqref{eq:relation_K}, also  $\widetilde{K}^{(N)}_{2,1}=\widetilde{K}^{(1)}_{2,1}$ for all $N\in\N$ and, finally, $J(N) =  NJ(1)$ thanks to the explicit formula~\eqref{eq:explicit_tildeK}.

To prove that $J(N)$ has no minimiser for $N\geq2$, we assume by contradiction that $\gamma$ is one.
By~\cite[Proposition~16]{GonLewNaz-20_ppt}, $\gamma$ is a rank $N$ projector. In addition, since we have equality in~\eqref{eq:relation_K}, $\gamma$ is also an optimiser for $K^{(N)}_{2,1}$. But then, by Theorem~\ref{thm:K}, it is of the form $\gamma=c\sum_{j=1}^N|\mu_j|^{1/2}\,|u_j\rangle\langle u_j|$ for some $c$. We conclude that $\mu_j=-1/c^2$ for all $j=1,...,N$ which is impossible since the first eigenvalue $\mu_1$ of a Schr\"odinger operator is always simple. 
\end{proof}

\begin{remark}\label{rem:ltconjj}
    In dimension $d=1$, a special case of the Lieb-Thirring conjecture~\cite{LieThi-76} states that 
    $$L_{\kappa,1}^{(N)}=L_{\kappa,1}^{(1)}\qquad\text{for all $\kappa\in(1,3/2]$ and all $N\geq1$.}$$
    If true, this conjecture would imply by the same argument as in the previous proof that 
    \begin{equation}
    J(N)=N\,J(1)\qquad \text{for all $2\leq p<3$ and all $N\geq1$, in dimension $d=1$,}
    \label{eq:conjecture_J_N_1}
    \end{equation}
and that the corresponding problems do not have minimisers for $N\geq 2$. The weaker conjecture~\eqref{eq:conjecture_J_N_1} appeared in~\cite{GonLewNaz-20_ppt}
\end{remark}

\subsection{Proof of Theorem~\ref{th:N=2}: triviality of solutions for $d = 1$, $p = 2$ and $N = 2$}
\label{ssec:proof_N=2}

In this subsection we prove Theorem~\ref{th:N=2}: we show that the fermionic NLS equation~\eqref{eq:fermionic_NLS} does not have a solution in the one dimensional case with $p = 2$ and $N = 2$. We will make use of the integrability of the equations. In the sequel, we study the ODE system 
 \begin{equation} \label{eq:NLS_N=2}
\begin{cases}
    v_1'' + 2 (v_1^2 + v_2^2) v_1 + \mu_1 v_1 = 0, \\
    v_2'' + 2 (v_1^2 + v_2^2) v_2 + \mu_2 v_2 = 0.
\end{cases}
\end{equation}
We added an extra factor $2$ to obtain the same explicit formulas as in the literature. If $(u_1, u_2)$ is a real-valued ground state solution to~\eqref{eq:no_binding_u}, then $(v_1, v_2) = \frac{1}{\sqrt{2}}(u_1, u_2)$ is a real-valued solution to~\eqref{eq:NLS_N=2}, which satisfies in addition $\| v_1 \| = \| v_2 \| = \frac12$.

The key step in the proof of Theorem~\ref{th:N=2} is the following classification result for~\eqref{eq:NLS_N=2} under an additional vanishing condition for $v_2$.

\begin{lemma} \label{lem:ODE}
    Let $\mu_1 \le \mu_2 < 0$, and let $(v_1, v_2)$ be a square integrable real-valued solutions of the ODE~\eqref{eq:NLS_N=2} with $v_2(0) = 0$. Then there are $a_1, a_2 \in \R$ such that
     \begin{equation}
\begin{cases}
          v_1(x) = \dfrac{a_1 \re^{ \eta_1x}}{f(x)}\left(1 +  \dfrac{ a_2^2}{4 \eta_2^2} \dfrac{\eta_1 - \eta_2}{\eta_1 + \eta_2}  \re^{2 \eta_2 x} \right),\\[0.4cm]
    v_2(x) = \dfrac{a_2 \re^{ \eta_2x}}{f(x)}\left(1 -  \dfrac{ a_1^2}{4 \eta_1^2} \dfrac{\eta_1 - \eta_2}{\eta_1 + \eta_2}  \re^{2 \eta_1 x} \right), 
\end{cases}
      \label{eq:explicit_1D}
     \end{equation}
    where
    \[
    f(x) = 1 + \dfrac{a_1^2}{4 \eta_1^2}  \re^{2 \eta_1 x} 
    + \dfrac{a_2^2}{4 \eta_2^2}  \re^{2 \eta_2 x}  +  \dfrac{a_1^2 a_2^2}{16 \eta_1^2 \eta_2^2} \dfrac{(\eta_1 - \eta_2)^2}{ (\eta_1 + \eta_2)^2} \re^{(2 \eta_2 + 2 \eta_1)x}
    \]
        and $\eta_1 := \sqrt{ | \mu_1 |}$, $\eta_2 := \sqrt{| \mu_2|}$.
\end{lemma}
In fact, if $a_2 \neq 0$, the condition $v_2(0) = 0$ fixes the value
\begin{equation} \label{eq:value_a1}
     a_1 = \pm 2 \eta_1 \left( \frac{\eta_1 + \eta_2}{\eta_1 - \eta_2}  \right)^{1/2}.
\end{equation}

\begin{proof}
    We proceed in two steps. First, we show that the functions~\eqref{eq:explicit_1D} are solutions and then we prove that they cover all possible initial data for $v_1(0)$, $v_1'(0)$ and $v_2'(0)$. By uniqueness of the solution of an initial value problem the result follows. 
    
    For the first point, checking the equation is simply a computation. For the convenience of the reader we quickly recall how to find the formulas~\eqref{eq:explicit_1D}. Following~\cite{RadLak-95} which uses Hirota's bilinearisation method~\cite{Hirota-80}, we write
    \[
    v_1 = \frac{g}{f}, \quad \text{and} \quad v_2 = \frac{h}{f}.
    \]
    With this change of variable, we see that~\eqref{eq:NLS_N=2} can we written as
    \[
    \begin{cases}
    f^2 \left( f g'' + f'' g - 2 f' g' + \mu_1 f g \right) + 2 f g \left( | f'|^2 -  f f''  + g^2 + h^2   \right) = 0, \\
    f^2 \left( f h'' + f'' h - 2 f' h' + \mu_2 f h \right) + 2 f h \left( | f'|^2 -  f f''  + g^2 + h^2   \right) = 0.
    \end{cases}
    \]
    We seek solutions that satisfy
    \[
    \begin{cases}
    f g'' + f'' g - 2 f' g' + \mu_1 f g = 0, \\
    f h'' + f'' h - 2 f' h' + \mu_2 f h = 0, \\
    | f'|^2 -  f f''  + g^2 + h^2 = 0.
    \end{cases}
    \]
    With Hirota's notation, this is of the form
    \[
    D(f,g) + \mu_1 fg = 0, \quad D(f,g) + \mu_2 fh = 0, \quad D(f, f) = \frac12 (g^2 + h^2),
    \]
    with the bilinear form $D(u,v) := uv'' + u''v - 2u'v'$. We now make the formal expansion $g = \chi g_1 + \chi^3 g_3$, $h = \chi h_1 + \chi^3 h_3$ and $f = 1 + \chi^2 f_2 + \chi^4$, and we solve the cascade of equations in powers of $\chi$. We first obtain (setting $\eta_1 := \sqrt{| \mu_1 |}$ and $\eta_2 := \sqrt{|\mu_2|}$)
    \[
    g_1 = a_1 \re^{ \eta_1 x}, \quad h_1 = a_2 \re^{ \eta_2 x}, 
    \]
    where $a_1$ and $a_2$ are two arbitrary constants. After some computation, we get (see also~\cite{RadLak-95}),
    \[
    f_2  = \dfrac{a_1^2}{4 \eta_1^2}  \re^{2 \eta_1 x} 
    + \dfrac{a_2^2}{4 \eta_2^2}  \re^{2 \eta_2 x},
    \]
    then
    \[
    g_3 = \left( \dfrac{a_1 a_2^2}{4 \eta_2^2} \dfrac{\eta_1 - \eta_2}{\eta_1 + \eta_2}   \right) \re^{(2 \eta_2 + \eta_1)x}, \quad
    h_3 = - \left( \dfrac{a_1^2 a_2}{4 \eta_1^2} \dfrac{\eta_1 - \eta_2}{\eta_1 + \eta_2}   \right) \re^{(2 \eta_1 + \eta_2)x}
    \]
    and finally
    \[
    f_4 = \dfrac{a_1^2 a_2^2}{16 \eta_1^2 \eta_2^2} \dfrac{(\eta_1 - \eta_2)^2}{ (\eta_1 + \eta_2)^2}\re^{(2 \eta_2 + 2 \eta_1)x}.
    \]
    This is the solution in Lemma~\ref{lem:ODE}. The condition $v_2(0) = 0$ gives the value of $a_1$ in~\eqref{eq:value_a1}.
    
    Let us now prove that all square integrable solutions with $v_2(0) = 0$ are of this form. In fact, instead of square integrability we will assume that $v_j$ and $v_j'$ tend to zero at infinity for $j=1,2$. It is not hard to deduce this property from the assumption that the solution is square integrable.
    
    For the proof we will assume that $v_2'(0)\neq 0$, for otherwise $v_2=0$ everywhere and the result is well-known (and easy to prove by a variation of the arguments that follow, using only \eqref{eq:cst_of_motion1} below).
       
    Any solution $(v_1, v_2)$ that decays at infinity has two constants of motion
    \begin{subequations} \label{eq:cst_of_motion}
        \begin{align}
        & (v_1^2 + v_2^2)^2  + | v_1' |^2 + | v_2 '|^2 + \mu_1 v_1^2 + \mu_2 v_2^2 = 0, \label{eq:cst_of_motion1} \\
        & (v_1^2 + v_2^2)(\mu_1 v_2^2 + \mu_2 v_1^2 + \mu_1 \mu_2) + (v_1 v_2' - v_1' v_2)^2 + \mu_2| v_1' |^2 + \mu_1 | v_2'|^2 = 0. \label{eq:cst_of_motion2}
        \end{align}
    \end{subequations}
	To obtain identity \eqref{eq:cst_of_motion1} we multiply the first and second equation in~\eqref{eq:NLS_N=2} by $v_1'$ and $v_2'$, respectively, add the resulting identities and then integrate using the fact that the solutions and their derivatives vanish at infinity. The fact that there is a second identity \eqref{eq:cst_of_motion2} reflects the integrability of the system~\cite{Manakov-74}.

 Evaluating~\eqref{eq:cst_of_motion} at $x = 0$ and using $v_2'(0) \neq 0$, we deduce that
    \[
        v_1(0)^2 = \mu_2 - \mu_1 \quad \text{and} \quad v_1'(0)^2 + v_2'(0)^2 = - \mu_2 \left( \mu_2 - \mu_1 \right).
    \]
    Thus, the value of $v_1(0)$ is determined, up to a sign, by $\mu_1$ and $\mu_2$ and we have
    \[
        v_1'(0)^2 < -\mu_2  (\mu_2 - \mu_1) = \eta_2^2 \left( \eta_1^2 - \eta_2^2 \right).
    \]
	The assumption $v_2'(0)\neq 0$ also shows that $- \mu_2(\mu_2-\mu_1) > 0$, hence $\mu_2 \neq \mu_1$ and therefore also $v_1(0)\neq 0$.

    Let $(\tilde v_1,\tilde v_2)$ be a solution of the form \eqref{eq:explicit_1D}. The absolute value of $a_1$ is fixed by~\eqref{eq:value_a1}. We will now show that the sign of $a_1$ as well as the number $a_2$ can be determined in such a way that $\tilde v_j(0)=v_j(0)$ and $\tilde v_j'(0)=v_j'(0)$ for $j=1,2$. Once we have shown this, ODE uniqueness implies that $\tilde v_j=v_j$ for $j=1,2$, which is what we wanted to prove. 
    
    Since $v_1(0)\neq 0$, we can choose the sign of $a_1$ in~\eqref{eq:value_a1} such that $\mathrm{sgn}\ a_1 = \mathrm{sgn}\ v_1(0)$. Note that, independently of the choice of $a_2$, we have $\mathrm{sgn}\ \tilde v_1(0)=\mathrm{sgn}\ a_1$. This, together with $\tilde v_1(0)^2 = \mu_2-\mu_1= v_1(0)^2$, implies that $\tilde v_1(0)= v_1(0)$.
    
    It remains to choose $a_2$. A tedious but straightforward computation yields
    \[
        \tilde v_1'(0) = - \frac{a_1}{|a_1|} \eta_2 \sqrt{\eta_1^2 - \eta_2^2}\ \dfrac{4 \eta_2^2(\eta_1 + \eta_2) - a_2^2  (\eta_1 - \eta_2)}{4 \eta_2^2(\eta_1 + \eta_2) + a_2^2(\eta_1 - \eta_2)}.
    \]
    The last quotient on the right side is a decreasing function of $a_2^2$ from $\left[ 0,  \infty \right]$ to $[-1,1]$. Thus, there is an $a_2^2\in (0,\infty)$ such that $\tilde v_1'(0)=v_1'(0)$. This determines the absolute value of $a_2$. To determine its sign, we note that the identities $\tilde v_1'(0)^2 + \tilde v_2'(0)^2 = - \mu_2 \left( \mu_2 - \mu_1 \right)= v_1'(0)^2 + v_2'(0)^2$ and $\tilde v_1'(0)=v_1'(0)$ imply that $\tilde v_2'(0)^2 = v_2'(0)^2$. Thus, we can choose the sign of $a_2$ in such a way that $\tilde v_2'(0) = v_2'(0)$.
    
    This shows that we can indeed find $a_1$ and $a_2$ such that $\tilde v_j(0)=v_j(0)$ and $\tilde v_j'(0)=v_j'(0)$ for $j=1,2$. As explained before, this implies the result.
\end{proof}

We will also need the following lemma in the proof of Theorem~\ref{th:N=2}.

\begin{lemma}\label{lem:normalization}
	If $(v_1, v_2)$ is a solution of the form~\eqref{eq:explicit_1D} of Lemma~\ref{lem:ODE}, then $\| v_1 \|^2 = 2 \eta_1$ and $\| v_2 \|^2 = 2 \eta_2$. In particular, we can have $\| v_1 \| = \| v_2 \|$ only if $\mu_1 = \mu_2$. 
\end{lemma}

\begin{proof}
	With the notation of Lemma~\ref{lem:ODE}, a computation reveals that
	\[
	v_1(x)^2 =   - \left( \frac{ \frac{a_2^2  \eta_1}{2 \eta_2^2} \re^{2 \eta_2 x} + 2 \eta_1 }{f(x)}    \right)'
	\quad \text{while} \quad
	v_2(x)^2 = - \left( \frac{\frac{a_1^2  \eta_2}{2 \eta_1^2} \re^{2 \eta_1 x} + 2 \eta_2 }{f(x)}   \right)'.
	\]
	Integrating gives
	\[
	\int_\R v_1^2 = -\left[  \frac{ \frac{a_2^2  \eta_1}{2 \eta_2^2} \re^{2 \eta_2 x} + 2 \eta_1 }{f(x)} \right]_{-\infty}^\infty = 2 \eta_1
	\quad \text{and similarly} \quad \int_\R v_2^2 = 2 \eta_2,
	\]
	as wanted.
\end{proof}

\begin{proof}[Proof of Theorem~\ref{th:N=2}]
	As explained before Lemma \ref{lem:ODE}, it is enough to consider solutions $(v_1,v_2)$ of~\eqref{eq:NLS_N=2} with $\| v_1 \| = \| v_2 \| = \frac12$.
	
	The equations~\eqref{eq:NLS_N=2} mean that the numbers $\mu_1$ and $\mu_2$ are negative eigenvalues of the operator $- \partial_{xx}^2 -2 (v_1^2+v_2^2)$. It is easy to see that the latter operator is bounded from below and its negative spectrum is discrete. Therefore it has a lowest eigenvalue $\mu_0$. Let $v_0$ be a corresponding eigenfunction, normalised by $\|v_0\|=\frac12$. It is well-known that the eigenvalue $\mu_0$ is non-degenerate and that $v_0$ can be chosen positive. In particular, if $v$ is a square integrable real valued solution to $- v'' -2 (v_1^2+v_2^2)v = \mu v$ which never vanishes, then necessarily $\mu=\mu_0$.
	
	We claim that $\mu_1=\mu_2=\mu_0$. To prove this, we may assume that $\mu_1\leq\mu_2<0$. In the case where $v_2$ never vanishes, the above remark gives $\mu_2=\mu_0$. Since $\mu_0$ is the lowest eigenvalue and since $\mu_1\leq\mu_2$, this also yields $\mu_1=\mu_0$. In the opposite case where $v_2$ does vanish at some point we can, after a translation, apply Lemma~\ref{lem:ODE}. We deduce that $v_1$ does not vanish, hence $\mu_1=\mu_0$. Moreover, applying Lemma~\ref{lem:normalization}, we conclude that $\mu_1=\mu_2$. This proves the claim.
	
	It follows from the equality $\mu_1=\mu_2=\mu_0$, the simplicity of $\mu_0$ and the normalisation that $v_1^2=v_2^2$. In particular, $v_1$ and $v_2$ both satisfy $v_j''+4v_j^3 + \mu_0 v_j = 0$. By uniqueness of the solution to the equation up to translations, this gives~\eqref{eq:formulas_u_1_u_2} for some $x_0\in\R$ and a sign $\pm$. Since $v_1^2=v_2^2$ the $x_0$'s for the two functions coincide, while the signs are independent. This completes the proof of the theorem.	
\end{proof}

\appendix

\section{Proof of Lemma~\ref{lem:duality_N}}
\label{appendix:proof_duality}

The proof of Lemma~\ref{lem:duality_N} splits naturally into two parts. We first deduce~\eqref{eq:LT_Schatten} from~\eqref{eq:LT_V_N}. We write our operator $\gamma$ in the form
$$\gamma=\sum_{j=1}^Nn_j|u_j\rangle\langle u_j|,  \quad \text{so that} \quad
\rho_\gamma(x) = \sum_{j=1}^N n_j | u_j |^2(x),
$$
where $(u_1,...,u_N)$ forms an orthonormal system. The inequality~\eqref{eq:LT_Schatten} which we wish to prove therefore reads 
	\begin{equation}\label{eq:lt_tildedual}
	\sum_{j=1}^Nn_j \|\nabla u_j\|^2
	\geq  K_{d,p}^{(N)}\left(\int_{\R^d} \rho_\gamma^p {\rd x} \right)^{\frac2{d(p-1)}} \left( \sum_{j=1}^N n_j^{q} \right)^{-\frac{2}{d(p-1)}+1}.
	\end{equation}
For a constant $\beta>0$ to be determined, let
	$$
	V(x) = - \beta \rho_\gamma^{p-1}.
	$$
	For $\kappa\geq 1$ we use H\"older's inequality in Schatten spaces~\cite{Simon-79} in the form
	$$
	\Tr AB \geq - \left( \sum_{n=1}^N \lambda_n(A)_-^\kappa \right)^{\frac1\kappa} \left( \Tr B^{\kappa'} \right)^{\frac{1}{\kappa'}}
	$$
	for all $B\geq 0$ of rank $\leq N$. Applying this with $A=-\Delta +V$ and $B=\gamma$ we obtain, in view of~\eqref{eq:LT_V_N},
	\begin{align*}
	& \sum_{j=1}^N n_j\int_{\R^d} |\nabla u_j|^2 \,\rd x - \beta \int_{\R^d} \left( \sum_{j=1}^N n_j|u_j|^2 \right)^p \rd x = \sum_{j=1}^Nn_j \int_{\R^d} \left( |\nabla u_j|^2 + V |u_j|^2\right)\rd x \\
	& \qquad \geq - \left( \sum_{j=1}^N n_j^{\kappa'}\right)^{\frac1{\kappa'}} \left( \sum_{j=1}^N \left| \lambda_j(-\Delta +V) \right|^\kappa\right)^{\frac1\kappa} \\
	& \qquad \geq - \norm{\gamma}_{\gS^{\kappa'}}\left( L_{\kappa,d}^{(N)} \int_{\R^d} V(x)_-^{\kappa+\frac{d}2}\,\rd x \right)^{\frac1\kappa} \\
	& \qquad = - \norm{\gamma}_{\gS^{\kappa'}} \left( L_{\kappa,d}^{(N)} \right)^{\frac1\kappa} \beta^{1+\frac{d}{2\kappa}} \left( \int_{\R^d} \rho_\gamma^{(p-1)\big(\kappa+\frac{d}2\big)}\rd x \right)^{\frac1\kappa}.
	\end{align*}
	We optimise in $\beta$ by choosing
	$$
	\beta = \left( \frac{2\kappa}{2\kappa+d} \, \frac{\int_{\R^d} \rho_\gamma^p \rd x}{\norm{\gamma}_{\gS^{\kappa'}} \left( L_{\kappa,d}^{(N)} \right)^{\frac1\kappa} \left( \int_{\R^d} \rho_\gamma^{(p-1)(\kappa+d/2)}\rd x \right)^{\frac1\kappa}} \right)^{\frac{2\kappa}d}
	$$
	and obtain
	\begin{equation*}
	\sum_{j=1}^N n_j \int_{\R^d} |\nabla u_j|^2 \,\rd x  \geq \left( \frac{2\kappa}{2\kappa+d} \right)^{\frac{2\kappa}d} \frac{d}{2\kappa+d} \frac{ \left( \int_{\R^d} \rho_\gamma^p\, \rd x \right)^{1+\frac{2\kappa}d}}{\|\gamma\|_{\gS^{\kappa'}}^{\frac{2\kappa}d} \left( L_{\kappa,d}^{(N)} \right)^{\frac2d} \left( \int_{\R^d} \rho_\gamma^{(p-1)(\kappa+\frac{d}2)}\rd x \right)^{\frac2d} } \,.
	\end{equation*}
	We now choose $\kappa=p' - d/2$, which is $> 1$ since $p<1+2/d$ and which ensures that $(p-1)(\kappa+d/2)=p$. Thus,
	$$
	\sum_{j=1}^N n_j\int_{\R^d} |\nabla u_j|^2 \,\rd x \geq \left( \frac{2p'-d}{2p'} \right)^{\frac{2p'}d - 1} \frac{d}{2p'} \, \frac{ \left( \int_{\R^d} \rho_\gamma^p\, \rd x \right)^{\frac{2}{d(p-1)}}}{\|\gamma\|_{\gS^{\kappa'}}^{\frac{2p'}d-1} \left( L_{p'-d/2,d}^{(N)} \right)^{\frac2d} } \,.
	$$
	Therefore, the best constant $K_{d,p}^{(N)}$ in \eqref{eq:lt_tildedual} satisfies
	$$
	K_{d,p}^{(N)} \geq \left( \frac{2p'-d}{2p'} \right)^{\frac{2p'}d - 1} \frac{d}{2p'} \, \frac{1}{\left( L_{p'-d/2,d}^{(N)} \right)^{\frac2d} } \,.
	$$
	
	Conversely, assume that inequality \eqref{eq:lt_tildedual} holds and let $V\in L^{\kappa+d/2}(\R^d)$. We assume that $-\Delta+V$ has at least $N$ negative eigenvalues, the other case being handled similarly. Let $u_1,\ldots,u_N$ be orthogonal eigenfunctions corresponding to the $N$ lowest eigenvalues of $-\Delta+V$ and let 
	$$\gamma=\sum_{j=1}^Nn_j\,|u_j\rangle\langle u_j|,\qquad 	n_j = |\lambda_j(-\Delta+V)|^{\kappa-1} \,,$$
	so that 
	$$
	\tr(-\Delta+V)\gamma = \sum_{j=1}^Nn_j\,\lambda_j(-\Delta+V)  = - \sum_{j=1}^N|\lambda_j(-\Delta+V)|^{\kappa} \,.
	$$
	We have, for $p$ such that $p' = \kappa + \frac{d}{2}$, 
	\begin{align*}
	\sum_{j=1}^N |\lambda_j(-\Delta+V)|^{\kappa} & = - \sum_{j=1}^N n_j\int_{\R^d} \left( |\nabla u_j|^2 + V |u_j|^2\right)\rd x \\
	& \leq -  K_{d,p}^{(N)} \left(\int_{\R^d} \rho_\gamma^p\,\rd x \right)^{\frac{2}{d(p-1)}} \left( \sum_{j=1}^N n_j^{\kappa'} \right)^{-\frac{2}{d(p-1)}+1}  + \left\| V_- \right\|_{p'} \left\| \rho_\gamma \right\|_p
	\end{align*}
    Setting $x := \| \rho \|_p$, this is of the form $ - \alpha x^{\frac{2p}{d(p-1)}} + \beta x$, with $\frac{2p}{d(p-1)} > 1$. So it is bounded from above by 
	\begin{align*}
	\sum_{j=1}^N |\lambda_j(-\Delta+V)|^{\kappa} & \leq \left( K_{d,p}^{(N)} \right)^{-\frac{d(p-1)}2({d+p-dp)}} \left( \frac{d}{2p'} \right)^{\frac{d}{2p'-d}} \left( \frac{2p'-d}{2p'} \right) \\
	& \quad \times \left( \int_{\R^d} V_-^{p'}\,\rd x \right)^{\frac{2}{2p'-d}} \left( \sum_{j=1}^N n_j^{\kappa'} \right)^{\frac{2-d(p-1)}{2p-d(p-1)}}.
	\end{align*}
    Recall that 
	$$
	n_j^{\kappa'} = |\lambda_n(-\Delta+V)|^{\kappa}
	$$
	and therefore the above inequality becomes
	$$
	\sum_{j=1}^N |\lambda_j(-\Delta+V)|^\kappa \leq \left( K_{d,p}^{(N)} \right)^{-\frac{d}{2}} \left( \frac{d}{2p'} \right)^{\frac{d}2} \left( \frac{2p'-d}{2p'} \right)^{\frac{2p'-d}2} \int_{\R^d} V_-^{p'}\,\rd x \,.
	$$
	Therefore the best constant $L_{\kappa,d}^{(N)}$ in \eqref{eq:LT_V_N} satisfies
	$$
	L_{\kappa,d}^{(N)} \leq \left( K_{d,p}^{(N)} \right)^{-\frac{d}{2}} \left( \frac{d}{2p'} \right)^{\frac{d}2} \left( \frac{2p'-d}{2p'} \right)^{\frac{2p'-d}2}.
	$$
	This proves the lemma.\qed

\section{Comments on the NLS model and its dual}
\label{app:NLS_comments}

This appendix contains some additional comments on the minimisation problem $J(\lambda)$ in~\eqref{eq:def_J} studied in~\cite{GonLewNaz-20_ppt}.
The latter was considered for $\lambda\in\R_+$ instead of just $\lambda=N\in\N$. It is equivalent to the inequality
\begin{multline}
\widetilde K_{p,d}^{(\lambda)} \left(\int_{\R^d}\rho_\gamma(x)^p\,\rd x\right)^{\frac{2}{d(p-1)}} \leq \Big(\tr(\gamma)\Big)^{\frac{d+2-dp}{d(p-1)}}\;\|\gamma\|^{\frac{2}{d}}\;\tr(-\Delta\gamma),\\ \text{for all $1\leq p\leq 1+\frac2d$}
 \label{eq:LT_sub_critical_all_lambda_bis}
\end{multline}
with $\tr(\gamma)\leq\lambda$, which is a particular interpolation between the trace formula $\| \gamma \|_{\gS^1} = \Tr(\gamma) =\| \rho_\gamma\|_1$, and the Lieb-Thirring inequality~\eqref{eq:K_forkappa=1} at $p=1+2/d$. As discussed in Subsection \ref{sec:duallt}, another interpolation involving the Schatten space norm $\|\gamma\|_q^{\frac{d+2-dp}{d(p-1)}+\frac{2}{d}}$ instead of $\norm{\gamma}_1^{\frac{d+2-dp}{d(p-1)}} \;\|\gamma\|^{\frac{2}{d}}$ is the dual Lieb-Thirring inequality~\eqref{eq:LT_Schatten_unconst}. 

\subsection{An inequality with no optimiser}
Optimising~\eqref{eq:LT_sub_critical_all_lambda_bis} over all possible $\lambda$'s, we arrive at the inequality without constraints
\begin{equation}
\boxed{\widetilde K_{p,d} \left(\int_{\R^d}\rho_\gamma(x)^p\,\rd x\right)^{\frac{2}{d(p-1)}} \leq \Big(\tr(\gamma)\Big)^{\frac{d+2-dp}{d(p-1)}}\;\|\gamma\|^{\frac{2}{d}}\;\tr(-\Delta\gamma),}
 \label{eq:LT_sub_critical_all_lambda}
\end{equation}
for all $\gamma=\gamma^*\geq0$, with the best constant
\begin{equation}
\widetilde{K}_{p,d}:=\left(\sup_{\lambda>0}\frac{|J(\lambda)|}{\lambda}\right)^{-\frac{d+2-dp}{d(p-1)}}\frac1{p-1}\left(\frac{d}{2p}\right)^{\frac{2}{d(p-1)}}\left(1+\frac2d-p\right)^{-\frac{d+2-dp}{d(p-1)}}.
\end{equation}

We recall from~\cite[Section~1.3]{GonLewNaz-20_ppt} that 
$$\sup_\lambda\frac{|J(\lambda)|}{\lambda}=\lim_{\lambda\to\ii}\frac{|J(\lambda)|}{\lambda}<\ii.$$
From the results in~\cite{GonLewNaz-20_ppt} we can deduce that the inequality~\eqref{eq:LT_sub_critical_all_lambda} has no optimiser.

\begin{lemma}\label{lem:no_optimiser_LT}
Let $d\geq1$ and $1<p<\min(2,1+2/d)$. Then $\widetilde K_{p,d}<\widetilde K_{p,d}^{(\lambda)}$ for all $\lambda>0$. In particular the inequality~\eqref{eq:LT_sub_critical_all_lambda} admits no optimiser. 
\end{lemma}

\begin{proof}
It was shown in~\cite[Corollary~22]{GonLewNaz-20_ppt} that $J(\lambda)/\lambda$ is always above its limit. Therefore $\widetilde K_{p,d}<\widetilde K_{p,d}^{(\lambda)}$ and there cannot be an optimiser with finite trace.
\end{proof}

We believe that the optimisers of $\widetilde K^{(N)}_{p,d}$ converge in the limit $N\to\ii$ to periodic or translation-invariant operators, as discussed at the end of Section~\ref{sec:LT} and in~\cite{GonLewNaz-20_ppt}.

\begin{remark}[Monotonicity in $p$]
By H\"older's inequality, for any $\gamma=\gamma^*\geq 0$ the function
$$
p\mapsto \left( \int_{\R^d} \rho_\gamma(x)^p \,dx \right)^{\frac 2{d(p-1)}} \left( \int_{\R^d} \rho_\gamma(x)\,dx \right)^{- \frac 2{d(p-1)}}
$$
is non-decreasing. This implies that $p\mapsto \widetilde K_{p,d}$ is non-increasing on the interval $(1,1+2/d)$. In particular, since $\widetilde K_{p,d}^{\rm sc}=K_{1+2/d,d}^{\rm sc}$ is independent of $p$, and $\widetilde K_{p,d} \ge K_{p,d}$, we deduce that if $\widetilde K_{p,d}=\widetilde K_{p,d}^{\rm sc}$ for some $p=p_0$, then the same inequality holds for all $1<p\leq p_0$. This generalises the observation in~\cite{GonLewNaz-20_ppt} that if the standard Lieb--Thirring conjecture holds for $\kappa=1$ (that is, $\widetilde K_{p,d}=\widetilde K_{p,d}^{\rm sc}$ for $p=1+2/d$), then $\widetilde K_{p,d}=\widetilde K_{p,d}^{\rm sc}$ for all $1<p<1+2/d$.
\end{remark}

\subsection{Dual inequality}

A natural question is to determine the inequality dual to~\eqref{eq:LT_sub_critical_all_lambda}. This is the object of the following lemma. 

\begin{lemma}[Dual formulation of~\eqref{eq:LT_sub_critical_all_lambda}]\label{lem:duality_NLS}
	Let $d\geq 1$ and let $\kappa>1$ and $p<1+2/d$ be related by $p'=\kappa+d/2$. Then~\eqref{eq:LT_sub_critical_all_lambda} is equivalent to the inequality
	\begin{eqnarray}
	\label{eq:lt_tilde}
	\boxed{\Tr (-\Delta+V+\tau)_- \leq \widetilde L_{\kappa,d}\; \tau^{1-\kappa} \int_{\R^d} V_-^{\kappa+\frac d2}\,\rd x, }
	\end{eqnarray}
    valid for all $\tau>0$ and all $V\in L^{\kappa+\frac d2}(\R^d)$, in the sense that the best constants are related by
    \begin{equation} \label{eq:duality_tildeK_and_tildeL}
    \widetilde  K_{p,d} \widetilde {L}_{\kappa, d}^{\frac{2}{d}} =  \left(1 - \frac{d(p-1)}{2}\right)^{\frac{d + 2 - dp}{d(p-1)}}   
    \frac{d}{2} \frac{(p-1)^{\frac{2 + d}{d}}}{p^{\frac{2p}{d(p-1)}}}
    = 
    \frac{d}{2}\dfrac{(\kappa - 1)^{\frac{2}{d}(\kappa - 1)}}{(\kappa + \frac{d}{2})^{\frac{2}{d}\kappa + 1}}.
    \end{equation}
\end{lemma}

\begin{proof}
	Assume that \eqref{eq:lt_tilde} holds and let $0\leq\gamma\leq 1$ of finite kinetic energy. Set $\lambda:=\Tr(\gamma)$ and $\rho:=\rho_\gamma$. Then, for all $\tau>0$ and all $0\geq V\in L^{\kappa+\frac d2}(\R^d)$, from \eqref{eq:lt_tilde} with the abbreviation $L:=\widetilde L_{p'-d/2,d}$ we have
	\begin{align*}
	\Tr(-\Delta\gamma) & = \Tr(-\Delta+V+\tau)\gamma - \int_{\R^d} V\rho\,\rd x - \tau\lambda \\
	& \geq - L \tau^{-\kappa+1} \int_{\R^d} V_-^{\kappa+\frac d2}\,\rd x + \int_{\R^d} V_- \rho\,\rd x - \tau\lambda \,.
	\end{align*}
	We first optimise in $V$ by taking
    \[
        V = - \frac{1}{L^{p-1}}  \frac{(p-1)^{p-1}}{p^{p-1}}  \tau^{(\kappa - 1)(p-1)} \rho^{p-1},
    \]
    and obtain
	$$
	\Tr(-\Delta\gamma) \geq \frac{(p-1)^{p-1}}{p^{p}} \frac{1}{L^{p-1}} \tau^{(\kappa-1)(p-1)} \int_{\R^d} \rho^p\,\rd x -\tau\lambda.
	$$
	We then optimise in $\tau$ by taking (note that $(\kappa - 1)(p-1) = 1 - \frac{d}{2}(p-1) \in (0, 1)$, so the function is indeed bounded from above)
    \[
        \tau = \frac{1}{\lambda^{\frac{2}{d(p-1)}}}  \left(1 - \frac{d(p-1)}{2} \right)^{\frac{2}{d(p-1)}} \frac{(p-1)^{\frac{2}{d}}}{p^{\frac{2p}{d(p-1)}}} \frac{1}{L^{\frac{2}{d}}} \left(\int_{\R^d} \rho^p\,\rd x \right)^{\frac{2}{d(p-1)}},
    \]
    and we obtain finally
    \[
        \Tr( - \Delta \gamma) \ge \frac{1}{\lambda^{\frac{d + 2 - dp}{d(p-1)}}}  \frac{1}{L^{\frac{2}{d}}}  \left(1 - \frac{d(p-1)}{2}\right)^{\frac{d + 2 - dp}{d(p-1)}}  
        \frac{d}{2} \frac{(p-1)^{\frac{2 + d}{d}}}{p^{\frac{2p}{d(p-1)}}}
         \left(\int_{\R^d} \rho^p\,\rd x \right)^{\frac{2}{d(p-1)}}.
    \]
    Comparing with~\eqref{eq:LT_sub_critical_all_lambda} shows the first bound
    \[
        \widetilde K_{p,d} L^{\frac{2}{d}} \geq  \left(1 - \frac{d(p-1)}{2}\right)^{\frac{d + 2 - dp}{d(p-1)}}   
        \frac{d}{2} \frac{(p-1)^{\frac{2 + d}{d}}}{p^{\frac{2p}{d(p-1)}}}.
    \]
	
	Conversely, assume that \eqref{eq:LT_sub_critical_all_lambda} holds and let $V\in L^{\kappa+\frac{d}{2}}(\R^d)$ and $\tau>0$. We set $\gamma=\1(-\Delta+V+\tau<0)$, $\rho=\rho_\gamma$ and $\lambda=\Tr(\gamma)$. We obtain, from \eqref{eq:LT_sub_critical_all_lambda} with the abbreviation $K=\tilde K_{p,d}$,
	\begin{align*}
	\Tr(-\Delta+V+\tau)_- & = - \Tr(-\Delta+V+\tau)\gamma = -\Tr(-\Delta\gamma) - \int_{\R^d} V\rho\, \rd x - \tau\lambda \\
	& \leq - K \frac{1}{\lambda^{\frac{d + 2 - dp}{d(p-1)}}} \left( \int_{\R^d} \rho^p\, \rd x \right)^{\frac{2}{d(p-1)}} + \int_{\R^d} V_-\rho\, \rd x - \tau\lambda \,.
	\end{align*}
	Seen as a function of $\lambda$, the right-hand side is smaller than its maximum, attained for
    \[
        \lambda = \left(\frac{2}{d(p-1)} - 1\right)^{\frac{d(p-1)}{2}} \left( \frac{K}{\tau} \right)^{\frac{d(p-1)}{2}}  \int_{\R^d} \rho^p \,\rd x \,,
    \]
    so
\begin{multline*}
	\Tr(-\Delta+V+\tau)_- \leq 	\int_{\R^d} V_-\rho\,\rd x\\ - \frac{2}{d(p-1)} \left( \frac{2}{d(p-1)} - 1 \right)^{\frac{d(p-1)}{2} - 1} K^{\frac{d(p-1)}{2}} \tau^{1 - \frac{d(p-1)}{2}  } \int_{\R^d} \rho^p \, \rd x .
\end{multline*}
	Now, seen as a function of $\rho$, it is again smaller than its maximum. We deduce that (recall that $\kappa = \frac{p}{p-1} - \frac{d}{2} = 1 + \frac{1}{p-1} + \frac{d}{2}$)
	\begin{multline*}
	\Tr(-\Delta+V+\tau)_-\\ \leq \left( \frac{d}{2} \right)^{\frac{1}{p-1}} \left( \frac{2}{d(p-1)} - 1 \right)^{\frac{d + 2 - dp}{2(p-1)}} \left( \frac {p-1}{p} \right)^{\frac p{p-1}}
	\frac{1}{K^{\frac d2}} \tau^{1 - \kappa} \int_{\R^d} V_-^{\kappa + \frac{d}{2}}\, \rd x \,.	 
	\end{multline*}
	Comparing with \eqref{eq:lt_tilde} shows  that
    \begin{align*}
        \widetilde L_{\kappa,d} K^{\frac d2} &  \leq  \left( \frac{d}{2} \right)^{\frac{1}{p-1}} \left( \frac{2}{d(p-1)} - 1 \right)^{\frac{d + 2 - dp}{2(p-1)}} \frac{(p-1)^{\frac{p}{p-1}}}{p^{\frac p{p-1}}} \\
        & =  \left( \frac{d}{2} \right)^{\frac{d}{2}} \left( 1 - \frac{d(p-1)}{2} \right)^{\frac{d + 2 - dp}{2(p-1)}} \frac{(p-1)^{1 + \frac{d}{2}}}{p^{\frac p{p-1}}} \,.
    \end{align*}
	This proves the lemma.
\end{proof}

\subsection{Weak Lieb-Thirring inequalities}
The dual inequality~\eqref{eq:lt_tilde} provides an estimate on the quantity
\begin{equation}
 \sup_{\tau>0}\left\{\tau^{\kappa-1}\Tr (-\Delta+V+\tau)_-\right\}=\sup_{\tau>0}\left\{\tau^{\kappa-1}\sum_{n\geq1}\Big(\lambda_n(-\Delta+V)+\tau\Big)_-\right\}.
 \label{eq:weak_lp_tau}
\end{equation}
A natural question is to ask how this supremum compares with 
$$\Tr (-\Delta+V)_-^\kappa=\sum_{n\geq1}|\lambda_n(-\Delta+V)|^\kappa$$
appearing in the usual Lieb-Thirring inequality. In this subsection we show that~\eqref{eq:weak_lp_tau} is equivalent to the weak $\ell^\kappa$ norm of the eigenvalues of $-\Delta+V$. In this sense~\eqref{eq:lt_tilde} is weaker than the ordinary Lieb-Thirring inequality for $\kappa$, which bounds the (strong) $\ell^\kappa$ norm of the eigenvalues. The results of this subsection concern the `analytic content' of the inequalities and ignore, at least to some extent, the question of sharp constants.

Let $X$ be a measure space and $p>r\geq 0$. For a measurable function $f$ we set
$$
[f]_{p,r}' := \sup_{\tau>0} \left\{\tau^{1-\frac rp} \left( \int_X (|f|-\tau)_+^r \,\rd x \right)^\frac 1p \right\}.
$$
When $r = 0$, we get
$$
[f]_{p,0}' = \sup_{\tau>0} \tau |\{|f|>\tau\}|^{1/p}
$$
which is the standard quasinorm in weak $L^p$. Actually, it turns out that for all $0 \le r<p$, $[f]_{p,r}'$ is an equivalent quasinorm in this space.

\begin{lemma}\label{weaknorms}
	If $p>r\geq 0$, then for any measurable $f$ on $X$,
	$$
	\left( \frac{(p-r)^{p-r}\, r^r}{p^p} \right)^{\frac 1p} [f]_{p,0}' \leq	
	[f]_{p,r}' \leq \left( \frac{\Gamma(p-r)\,\Gamma(r+1)}{\Gamma(p)} \right)^\frac 1p [f]_{p,0}' \,.
	$$
\end{lemma}

\begin{proof}
	We set $\lambda(\sigma) := |\{|f|>\sigma\} |$ for brevity. First, for any $\sigma>\tau$, we have the inequality
    \[
        \1_{ \{| f | > \sigma \} } \le \1_{\{| f | > \sigma\} } \left( \dfrac{| f | - \tau}{\sigma - \tau} \right)^r 
        \le  \1_{\{| f | > \tau\} } \left( \dfrac{| f | - \tau}{\sigma - \tau} \right)^r.
    \]
    Integrating gives the inequality
	$$
	\lambda(\sigma) \leq \dfrac{1}{(\sigma-\tau)^{r}} \int_X (|f|>\tau)_+^r\,\rd x \leq \dfrac{1}{\tau^{p-r} (\sigma-\tau)^{r}} \left( [f]_{p,r}' \right)^p \,.
	$$
	We optimise in $\tau$ by choosing $\tau = \left(\frac{p-r}{p} \right) \sigma$, and obtain that
    \begin{equation*}
        \sigma^p \lambda(\sigma) \le \dfrac{p^p}{(p-r)^{p-r} r^r}  \left( [f]_{p,r}' \right)^p.
    \end{equation*}
	which is the first bound. Conversely, we use the identity
	$$
	(|f|-\tau)_+^r = r \int_\tau^\infty \1_{\{ |f|>\sigma\}} (\sigma-\tau)^{r-1} \,\rd\sigma.
	$$
	Integrating over $X$ gives
	\begin{eqnarray}
	\label{eq:layercake}
		\tau^{p-r} \int_X (|f|-\tau)_+^r \,\rd x = r \tau^{p-r} \int_\tau^\infty \lambda(\sigma) (\sigma-\tau)^{r-1} \,\rd\sigma \,.
	\end{eqnarray}
	Estimating $\lambda(\sigma)\leq \sigma^{-p} \left( [f]_{p,0}' \right)^p$ we obtain
	$$
	\tau^{p-r} \int_X (|f|-\tau)_+^r \,\rd x \leq r \left( [f]_{p,0}' \right)^p \int_1^\infty \dfrac{(s-1)^{r-1}}{s^p}\,\rd s =  \left( [f]_{p,0}' \right)^p \frac{r\,\Gamma(p-r)\,\Gamma(r)}{\Gamma(p)} \,,
	$$
	which is the second bound.
\end{proof}

Note that if $\lambda_n(-\Delta+V)$ denote the negative eigenvalues of $-\Delta+V$, repeated according to multiplicities, then
$$
\sup_{\tau>0} \tau^{\kappa-1} \Tr(-\Delta+V+\tau)^\kappa_- = \left( \left[ \lambda_\cdot(-\Delta+V) \right]_{\kappa,1}' \right)^{\frac 1\kappa} \,.
$$
Thus, combining Lemmas \ref{lem:duality_NLS} and \ref{weaknorms}, we obtain

\begin{corollary}[Weak Lieb-Thirring inequality]
	Inequalities \eqref{eq:lt_tilde} and \eqref{eq:LT_sub_critical_all_lambda} are equivalent to the inequality
	$$
	\boxed{\norm{ \Big(\lambda_n (-\Delta+V)\Big)_{n\geq1} }_{\ell^\kappa_{\rm w}}^\kappa \lesssim \int_{\R^d} V(x)_-^{\kappa+\frac d2}\, \rd x  }
	$$
	for all $V\in L^{\kappa+\frac d2}(\R^d)$.
\end{corollary}

The equivalence claimed in this corollary is weaker than that in Lemma \ref{lem:duality_NLS} since the (not displayed) constant depends on the choice of the norm in $\ell^\kappa_{\rm w}$.

\subsection{Semiclassical constants}

It was proved in~\cite[Lemma~10]{GonLewNaz-20_ppt} that $\widetilde K_{p,d}$ is not larger than its semiclassical counterpart, which is independent of $p$ and given by the $p=1+2/d$ semi-classical constant
$$
\widetilde K_{d}^{\rm sc}= K_{1+2/d,d}^{\rm sc} = \frac{4 \pi^2 d}{d+2} \left( \frac{d}{| \SS^{d-1} |} \right)^{\frac{2}{d}}.
$$
Together with Proposition~\ref{prop:relation_K}, we obtain
$$\boxed{ K_{p,d} \le \widetilde K_{p,d}\leq \widetilde K_{d}^{\rm sc} .} $$
In the dual picture, we have a similar result:
\begin{lemma} \label{lem:sc_L}
    For all $\kappa \ge 1$, we have
    \begin{equation} \label{eq:ineq_for_L_tildeL}
        \boxed{ \dfrac{(\kappa - 1)^{\kappa - 1}}{\kappa^\kappa} \ L_{\kappa, d} \ge \widetilde L_{\kappa,d}\geq \widetilde L_{\kappa, d}^{\rm sc}, }
    \end{equation}
    where the semi-classical constant $\tilde L_{\kappa, d}^{\rm sc}$ is defined by
    \begin{equation} \label{eq:explicit_tildeLsc}
        \widetilde L_{\kappa, d}^{\rm sc} := \dfrac{(\kappa - 1)^{(\kappa - 1)}(1 + \frac{d}{2})^{1 + \frac{d}{2}}}{(\kappa + \frac{d}{2})^{\kappa + \frac{d}{2}}}  L^{\rm sc}_{1,d}
    \end{equation}
    with the semiclassical constant $L^{\rm sc}_{1,d}$ at $\kappa=1$ given by~\eqref{eq:L_sc}.
\end{lemma}

\begin{proof}
	Both inequalities in~\eqref{eq:ineq_for_L_tildeL} follow from the explicit formulas~\eqref{eq:duality_K_and_L_unconst} and~\eqref{eq:duality_tildeK_and_tildeL}.
\end{proof}

\begin{remark}[The semi-classical constant]
	We show here that the constant $\widetilde L_{\kappa, d}^{\rm sc}$ has an interpretation in terms of a semiclassical limit, thereby justifying its name. Because of the second inequality in \eqref{eq:explicit_tildeLsc}, this argument shows that considered scenarios is in a certain sense dual to that considered in \cite{GonLewNaz-20_ppt}. For any $V\in L^{\kappa+\frac d2}(\R^d)$ and any $\tau>0$, we have
	$$
	\tau^{\kappa-1} \Tr(-\hbar^2\Delta+V+\tau)_- \underset{\hbar\to0}\sim \tau^{\kappa-1} \hbar^{-d} L_{1,d}^{\rm sc} \int_{\R^d} (V+\tau)_-^{1+\frac d2}\,\rd x \,.
	$$
	On the other hand, by inequality \eqref{eq:lt_tilde},
	\begin{align*}
	\tau^{\kappa-1} \Tr(-\hbar^2\Delta+V+\tau)_- & = \hbar^{2\kappa} ( \hbar^{-2}\tau)^{\kappa-1} \Tr(-\Delta+\hbar^{-2} V + \hbar^{-2} \tau)_- \\
	& \leq \hbar^{2\kappa} \widetilde L_{\kappa,d} \int_{\R^d} \left( \hbar^{-2} V \right)_-^{\kappa+\frac d2}\,\rd x = \hbar^{-d} \widetilde L_{\kappa,d} \int_{\R^d} V_-^{\kappa+\frac d2}\,\rd x \,.
	\end{align*}
	This shows that
	$$
	\tau^{\kappa-1}  \int_{\R^d} (V+\tau)_-^{1+\frac d2}\,dx
	\leq \frac{\widetilde L_{\kappa,d}}{L_{1,d}^{\rm sc}} \int_{\R^d} V_-^{\kappa+\frac d2}\,\rd x \,.
	$$
	Taking the supremum in $\tau$ shows that
	\[
	\left[V_- \right]'_{\kappa+\frac d2,1+\frac d2}  \le \left( \frac{\widetilde L_{\kappa,d}}{L_{1,d}^{\rm sc}} \right)^{\frac{1}{\kappa+\frac d2}} \left\| V_- \right\|_{L^{\kappa+\frac d2}}.
	\]
	According to the optimality statement in the following lemma, we have
	$$
	\left( \frac{\widetilde L_{\kappa,d}}{L_{1,d}^{\rm sc}} \right)^{\frac{1}{\kappa+\frac d2}} \geq \left( \frac{(\kappa-1)^{\kappa-1}\, (1+\frac d2)^{1+\frac d2}}{(\kappa+\frac d2)^{\kappa+\frac d2}} \right)^\frac 1{\kappa+\frac d2}.
	$$
	This proves, once again, the second inequality in \eqref{eq:explicit_tildeLsc} and shows how this inequality is related to a semiclassical limit.
\end{remark}	
	
	\begin{lemma}\label{lem:weakstrong}
		Let $X$ be a measure space, $p>r\geq 0$ and $f\in L^p(X)$. Then
		$$
		[f]_{p,r}' \leq \left( \frac{(p-r)^{p-r}\,r^r}{p^p} \right)^\frac 1p \|f\|_p \,.
		$$
		The constant on the right side is best possible.
	\end{lemma}
	
	\begin{proof}
		We first recall that
		$$
		\int_X |f|^p\,\rd x = p \int_0^\infty \lambda(\sigma)\sigma^{p-1}\,\rd\sigma.
		$$
		Together with~\eqref{eq:layercake} (note that we may assume $r > 0$ by continuity) we need to prove that
		$$
		r \tau^{p-r} \int_\tau^\infty \lambda(\sigma) (\sigma-\tau)^{r-1} \,\rd\sigma
		\leq \frac{(p-r)^{p-r}\,r^r}{p^p} p \int_0^\infty \lambda(\sigma)\sigma^{p-1}\,\rd\sigma \,.
		$$
		We write $\lambda = \int_0^\infty \1_{\{\lambda> b\}}\,\rd b$ and, since $\lambda$ is non-increasing, for any $b>0$ the function $\1_{\{\lambda> b\}}$ is the characteristic function of an interval with left endpoint at zero. Thus, it suffices to prove the above inequality for such characteristic functions. A computation shows that
		$$
		r \tau^{p-r} \int_\tau^\infty \1_{[0,a)}(\sigma)(\sigma-\tau)^{r-1} \,\rd\sigma = \tau^{p-r} (a-\tau)_+^r$$
		and
		$$	p \int_0^\infty \1_{[0,a)}(\sigma)\sigma^{p-1}\,\rd\sigma = a^p \,.
		$$
		Thus, the inequality follows from the elementary equality
		$$
		\sup_{a>0} \tau^{p-r} (a-\tau)_+^r  = \frac{(p-r)^{p-r}\,r^r}{p^p} a^p \,.
		$$	
		There is equality when $f$ is a characteristic function and $\tau$ is chosen appropriately. This proves Lemma~\ref{lem:weakstrong}.
	\end{proof}

\begin{remark}
    We wonder whether for all $d \ge 1$ and all $\kappa \ge \frac32$, we have the equality $\widetilde L_{\kappa, d} = \widetilde L_{\kappa, d}^{\rm sc}$. This would be the analogue of the equality $L_{\kappa, d} = L_{\kappa, d}^{\rm sc}$~\cite{LapWei-00}. We have the following rather tight bounds. Thanks to the explicit formulas~\eqref{eq:explicit_tildeLsc} and~\eqref{eq:L_sc}, one can numerically plot the two curves $\kappa \mapsto \widetilde L_{\kappa, d}^{\rm sc}$ and $\kappa \mapsto \dfrac{(\kappa - 1)^{\kappa - 1}}{\kappa^\kappa} L_{\kappa, d}^{\rm sc}$. As stated in Lemma~\ref{lem:sc_L}, the two curves coincide at $\kappa = 1$, but for all $\kappa > 1$, it appears that
    \[
        0 < \dfrac{(\kappa - 1)^{\kappa - 1}}{\kappa^\kappa}\, L_{\kappa, d}^{\rm sc} - \widetilde L_{\kappa, d}^{\rm sc} < \begin{cases}
            0.004 & \text{for} \quad d = 1, \\
            0.0009 & \text{for} \quad d = 2, \\
            0.0002 & \text{for} \quad d = 3. \\
        \end{cases}
    \]
    In the region $\kappa \ge 3/2$ where $L_{\kappa, d} = L_{\kappa, d}^{\rm sc}$~\cite{LapWei-00}, we deduce that $| \widetilde L_{\kappa, d} - \widetilde L_{\kappa, d}^{\rm sc}|$ is smaller than the constants above.
\end{remark}

\section{An inequality on the other side of the Lieb-Thirring exponent}
\label{app:LT-Sobolev}

In this section we would like to compare our inequality~\eqref{eq:LT_sub_critical_all_lambda} with the following related inequality,
\begin{multline}
K'_{p,d} \norm{\rho_\gamma}_{L^{p}(\R^d)}^\frac{2p}{d(p-1)}\leq \|\gamma\|^{\frac{d-(d-2)p}{d(p-1)}}\;\tr(-\Delta\gamma),\\ 1+\frac2d\leq p< 1+\frac{2}{d-2},\quad d\geq3.
\label{eq:LT_super_critical}
\end{multline}
This inequality remains valid in dimensions $d=1,2$, with $1/(d-2)$ replaced by $+\ii$. Note that the exponent $p$ in \eqref{eq:LT_super_critical} lies on the other side of the Lieb-Thirring exponent, compared to the situation considered in this paper. Inequality \eqref{eq:LT_super_critical} appears in~\cite[Eq.~(3.7)]{LieLos-86} for $p=2$ and $d=3$. 

The proof of \eqref{eq:LT_super_critical} in dimension $d\geq 3$ is simple. Indeed, the Hoffmann-Ostenhof~\cite{Hof-77} inequality~\eqref{eq:Hoffmann-Ostenhof} together with the Sobolev inequality give
\begin{equation}
S_{\frac{d}{d-2},d}\norm{\rho_\gamma}_{L^{\frac{d}{d-2}}(\R^d)}\leq \tr(-\Delta\gamma)\qquad\text{for all $d\geq3$ and all $\gamma=\gamma^*\geq0$.}
 \label{eq:Sobolev}
\end{equation}
Using H\"older's inequality and the Lieb-Thirring inequality~\eqref{eq:K_forkappa=1} (with constant $K_{1+2/d,d}=\inf_N K_{1+2/d,d}^{(N)}>0$) we obtain \eqref{eq:LT_super_critical}.

Our inequality~\eqref{eq:LT_sub_critical_all_lambda} interpolates between the Lieb-Thirring inequality and the trace equality $\| \gamma \|_{\gS^1} = \Tr(\gamma) =\| \rho_\gamma\|_1$. In contrast, inequality~\eqref{eq:LT_super_critical} interpolates between the Lieb-Thirring inequality \eqref{eq:K_forkappa=1} and the Sobolev inequality~\eqref{eq:Sobolev}. 

An interesting difference between \eqref{eq:LT_sub_critical_all_lambda} and \eqref{eq:LT_super_critical} arises when one considers the question of existence of minimizers. Recall from Lemma~\ref{lem:no_optimiser_LT} that~\eqref{eq:LT_sub_critical_all_lambda} never has optimisers. On the other hand, in~\cite{HonKwoYoo-19} the existence of optimisers for~\eqref{eq:LT_super_critical} was proved when $1+2/d<p< 1+2/(d-2)$. When normalised in the manner $\|\gamma\|=1$ and $\tr(-\Delta\gamma)=\theta\int_{\R^d}\rho_\gamma^p$, these optimisers were shown in~\cite[Thm.~2]{HonKwoYoo-19} to solve the equation
\begin{equation}
\gamma=\1_{(-\ii,0)}\left(-\Delta-\rho_\gamma^{p-1}\right)+\delta,\qquad\text{with}\quad  0\leq\delta=\delta^*\leq\1_{\{0\}}\left(-\Delta-\rho_\gamma^{p-1}\right). 
\label{eq:HKY_equation}
\end{equation}
In other words, $\gamma$ is the orthogonal projection onto all the negative eigenfunctions, except possibly on the kernel of $-\Delta-\rho_\gamma^{p-1}$. If these optimisers $\gamma$ have a finite rank $N$ (they do for $d\geq3$ and $p$ large enough), then they must be NLS ground states in the sense of~\cite{GonLewNaz-20_ppt}.

We now slightly refine the result in~\cite{HonKwoYoo-19} by showing that the operator $-\Delta-\rho_\gamma^{p-1}$ has no zero eigenvalues and, in particular, one has $\delta=0$ in \eqref{eq:HKY_equation}.

\begin{proposition}
	Let $1+\frac{2}{d}< p<\infty$ if $d=1,2$ and $1+\frac{2}{d}<p<1+\frac{2}{d-2}$ and let $\gamma$ be an optimiser for~\eqref{eq:LT_super_critical}, normalised so that $\|\gamma\|=1$ and $\tr(-\Delta\gamma)=\theta\int_{\R^d}\rho_\gamma^p$. Then
	$$
	\ker\left(-\Delta-\rho_\gamma^{p-1}\right)=\{0\} \,.
	$$
\end{proposition}

\begin{proof}
	We begin by proving that $\delta=0$ in \eqref{eq:HKY_equation}. We denote by $u_j$ and $\mu_j$ the eigenfunctions and eigenvalues of $-\Delta-\rho_\gamma^{p-1}$ and by $n_j$ the corresponding eigenvalues of $\gamma$. From~\eqref{eq:HKY_equation} we know that $n_j=1$ if $\mu_j<0$. By arguing as in~\eqref{eq:estim_mu_j}, we have the estimate
	\begin{equation}
	\mu_{j}\leq \frac{\theta\int_{\R^d}\rho_\gamma^p}{n_j}\left(1-\frac{n_j}{\theta\int_{\R^d}\rho_\gamma^p}\int_{\R^d}\rho_\gamma^{p-1}|u_j|^2-\left(\frac{\int_{\R^d}(\rho_\gamma-n_j|u_j|^2)^p}{\int_{\R^d}\rho_\gamma^p}\right)^{\frac1{\theta p}}\right)
	\label{eq:estim_mu_j_HKY}
	\end{equation}
	with $\theta=d/(2p')\in(1/p,1)$. We claim that the right side is negative, which yields $\mu_j<0$, that is, $\delta\equiv0$ in~\eqref{eq:HKY_equation}. To see this, we remark that for any $f\geq0$ and any probability measure $\mathbb{P}$, we have by H\"older's inequality twice
	$$\int f\,d\mathbb{P}\leq \left(\int f^p\,d\mathbb{P}\right) ^{\frac1p}\leq \theta \left(\int f^p\,d\mathbb{P}\right) ^{\frac1{\theta p}}+(1-\theta).$$
	The second inequality is strict when $\int f^p\,d\bP\neq 1$. This may be rewritten in the form
	\begin{equation}
	1+\theta^{-1}\int (f-1)\,d\mathbb{P}\leq \left(\int f^p\,d\mathbb{P}\right) ^{\frac1{\theta p}}. 
	\label{eq:Holder_f_P}
	\end{equation}
	Choosing $f=1-n_j|u_{j}|^2/\rho_\gamma$ and $\mathbb{P}=\rho_\gamma^p/\int_{\R^d}\rho_\gamma^p$, we obtain $\mu_j<0$ in~\eqref{eq:estim_mu_j_HKY} since $f\leq1$ and $f\neq1$, hence $\int_{\R^d}f^p\,d\bP<1$.
	We have thus proved that  $\delta\equiv0$ in~\eqref{eq:HKY_equation}.
	
	We now show that $\ker(-\Delta-\rho_\gamma^{p-1})=\{0\}$. Indeed, assume on the contrary that $\mu_j=0$ (then $n_j=0$ by the previous argument). Consider this time the perturbation $\gamma(t)=\gamma+t|u_j\rangle\langle u_j|$, which cannot be an optimiser for $t>0$. Taking $\mu_j = 0$ and $n_j = -t$ in~\eqref{eq:estim_mu_j_HKY} gives the (strict) inequality
	\begin{equation}
	\left(\frac{ \int_{\R^d}\left(\rho_\gamma+t|u_j|^2\right)^p}{ \int_{\R^d}\rho_\gamma^p}\right)^{\frac1{\theta p}}<1+\frac{t\int_{\R^d}\rho_\gamma^{p-1}|u_j|^2}{\theta\int_{\R^d}\rho_\gamma^p}
	\label{eq:estim_mu_j_HKY_zero}
	\end{equation}
	for all $0< t<\|\gamma\|$. By~\eqref{eq:Holder_f_P} with $f=1+t|u_j|^2/\rho_\gamma$, which satisfies $\int_{\R^d}f^p\,d\bP>1$, we have
	$$\left(\frac{ \int_{\R^d}\left(\rho_\gamma+t|u_j|^2\right)^p}{ \int_{\R^d}\rho_\gamma^p}\right)^{\frac1{\theta p}}>1+\frac{t\int_{\R^d}\rho_\gamma^{p-1}|u_j|^2}{\theta\int_{\R^d}\rho_\gamma^p}$$
	and we obtain a contradiction. Therefore $\ker(-\Delta-\rho_\gamma^{p-1})=\{0\}$, as claimed.
\end{proof}

\begin{figure}[h]
\includegraphics[width=10cm]{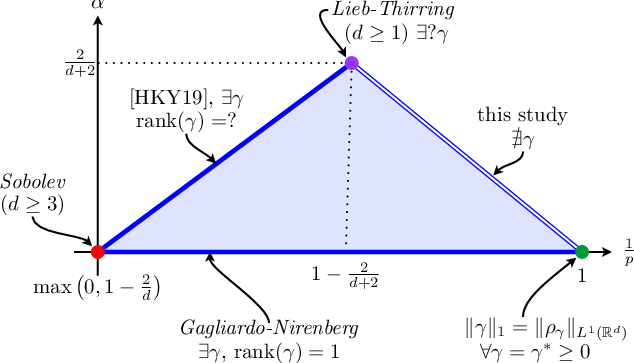}

\caption{Graphical representation of the validity and existence of optimisers for Lieb-Thirring-type inequalities in the form 
$$\norm{\rho_\gamma}_{L^p(\R^d)}\leq C\|\gamma\|^\alpha\norm{\gamma}_1^{\beta}\norm{\sqrt{-\Delta}\gamma\sqrt{-\Delta}}_1^{1-\alpha-\beta}.$$ 
We deal in~\cite{GonLewNaz-20_ppt} and this paper with the right edge where $\alpha,\beta>0$. There is no optimiser without an additional trace constraint. Existence of optimisers was proved on the left edge where $\beta=0$ in~\cite{HonKwoYoo-19}. The horizontal edge coincides with the Gagliardo-Nirenberg inequality, with $\alpha=0$. Minimisers exist and are all rank-one. In dimension $d\geq3$, the Sobolev inequality has a formal rank-one optimiser. For $d=3,4$, however, it is not bounded on $L^2(\R^d)$ since the associated function is not in $L^2(\R^d)$.  It is expected that a minimiser exists for the Lieb-Thirring inequality only in dimension $d=1$, where it should be rank-one. In dimension $d=1$, our study is limited to $p<2$. 
\label{fig:inequalities}}
\end{figure}


\subsection*{Acknowledgement} This project has received funding from the U.S. National Science Foundation (grant agreements DMS-1363432 and DMS-1954995 of R.L.F.) and from the European Research Council (ERC) under the European Union's Horizon 2020 research and innovation programme (grant agreement MDFT 725528 of M.L.).


\end{document}